\documentclass[twocolumn]{autart}    

\usepackage{amsmath,amsfonts}
\usepackage{algorithm}
\usepackage{algpseudocode}
\usepackage{array}
\usepackage{caption}
\usepackage{textcomp}
\usepackage{stfloats}
\usepackage{subcaption}
\usepackage{graphicx}
\usepackage{upgreek}
\usepackage{url}
\usepackage{verbatim}
\usepackage{graphicx}
\usepackage{mathtools}
\usepackage{arydshln}
\usepackage{ulem}
\usepackage{cite}
\usepackage{xcolor}

\begin{document}

\begin{frontmatter}

\title{Finite Control Set Model Predictive Control\\ with Limit Cycle Stability Guarantees\thanksref{footnoteinfo}} 

\thanks[footnoteinfo]{This work is funded by the EU Horizon 2020 research project IT2 (IC Technology for the 2nm Node), Grant agreement ID: 875999. }

\author[Eind]{Duo Xu}\ead{d.xu@tue.nl},    
\author[Eind]{Mircea Lazar}\ead{m.lazar@tue.nl}               

\address[Eind]{Eindhoven University of Technology, Eindhoven}  


\begin{abstract}                          
This paper considers the design of finite control set model predictive control (FCS-MPC) for discrete-time switched affine systems. Existing FCS-MPC methods typically pursue practical stability guarantees, which ensure convergence to a bounded invariant set that contains a desired steady state. As such, current FCS-MPC methods result in unpredictable steady-state behavior due to arbitrary switching among the available finite control inputs. Motivated by this, we present a FCS-MPC design that aims to stabilize a steady-state limit cycle compatible with a desired output reference via a suitable cost function. We provide conditions in terms of periodic terminal costs and finite control set control laws that guarantee asymptotic stability of the developed limit cycle FCS-MPC algorithm. Moreover, we develop conditions for recursive feasibility of limit cycle FCS-MPC in terms of periodic terminal sets and we provide systematic methods for computing ellipsoidal and polytopic periodically invariant sets that contain a desired steady-state limit cycle. Compared to existing periodic terminal ingredients for tracking MPC with a continuous control set, we design and compute terminal ingredients using a finite control set. The developed methodology is validated on switched systems and power electronics benchmark examples.

\end{abstract}

\end{frontmatter}

\section{Introduction}

Over past decades, switched systems have witnessed remarkable developments from the theoretical point of view to practical insights, thereby offering substantial benefits in numerous real-world applications \cite{benmiloud2018local}. Among switched systems applications, power electronics play a crucial role in many domains, such as power conversion in electrical vehicles charging and renewable energy systems \cite{deaecto2016stability}.
More specifically, this paper focuses on a specific category of switched systems that comprises a finite number of subsystems, each characterized by a constant state matrix and constant affine term. These subsystems are governed by a switching rule that regulate the transitions between them \cite{liberzon2003switching}. 

The presence of the affine terms create several equilibrium points corresponding to each subsystem, which poses a notable challenge in designing effective switching control laws aimed at stabilization. This difficulty primarily stems from the fact that the desired reference point for stabilization often does not coincide with the equilibrium point of any isolated subsystem \cite{sun2011stability}. In general, for switched affine systems, or equivalently, for finite control set bilinear affine systems, the attainable reference point can be determined using an averaged model. This model represents the aggregated performance across a designated set of subsystems. 
Consequently, the practical stability approach for a fixed attainable reference point has been extensively explored within the context of switched affine systems, leading to significant advancements in the continuous-time setting \cite{serieye2023attractors}.
In the discrete-time setting, \cite{deaecto2016stability} focuses on the practical stability analysis and control design of switched affine systems using time-invariant Lyapunov functions.

In addition, several finite control set model predictive control (FCS-MPC) related papers have addressed this concern by providing practical stability (ultimate boundedness) conditions for a specific class of switched affine systems, which can be modeled as linear time-invariant systems with quantized inputs \cite{aguilera2011finite,aguilera2011stability}.
Therein, terminal ingredients for FCS-MPC were designed based on \textit{continuous control set} local control laws and the effect of using quantization functions to obtain quantized inputs was modeled as an additive disturbance. This facilitated practical stability conditions \cite{aguilera2013stability,quevedo2019model} that ensure convergence of the closed-loop FCS-MPC trajectories to a bounded invariant set compatible with a desired state reference. 

However, in many applications such as power electronics, the guarantees provided by practical stability are often insufficient to meet the required steady-state specifications. These specifications typically require a desired average switching frequency and total harmonic distortion, which cannot be guaranteed by practical stability results. Heuristic solutions meant to cope with such specifications via FCS-MPC include penalizing the output tracking error and incorporating penalties for the input rate of change \cite{prior2015control}.
In the study outlined in \cite{karamanakos2019guidelines}, a comparison was made between the output error formulations based on the 1-norm and the 2-norm. Additionally, for three-phase applications, \cite{geyer2018algebraic} demonstrates that penalizing the 2-norm of current errors in alpha-beta coordinates is directly equivalent to penalizing the total harmonic distortion.

Motivated by power electronics specifications and open problems in stabilization of switched affine systems, recent works have considered the stabilization of an adequately designed closed trajectory for switched affine systems, namely a limit cycle, instead of practical stabilization of a fixed steady state point. Conditions for asymptotic stability of a limit cycle with fixed periodicity are proposed together with a time-varying Lyapunov function in \cite{egidio2020global}. Based on these conditions, a time-dependent min-switching control rule is determined therein to regulate the states towards the limit cycle. Alternatively, a pure state-feedback min-switching control strategy is developed in \cite{serieye2023attractors}, together with robust stabilization criteria for uncertain switched affine systems. Inspired by these recent results on limit cycle stabilization, we developed a preliminary limit cycle tracking FCS-MPC formulation in \cite{xu2022steady}. Therein, we considered the specific class of linear systems with a finite control set and we established conditions for limit cycle asymptotic convergence in terms of terminal cost design, without any recursive feasibility guarantees.

In this paper we present a new approach to the design of FCS-MPC that aims to stabilize a desired limit cycle and constructs stabilizing terminal ingredients in the presence of a finite control set. In summary, the main contributions with respect to existing FCS-MPC frameworks, including our preliminary work \cite{xu2022steady}, are: 
\begin{itemize}
\item[\textit{(i)}] Design of FCS-MPC for limit cycle stabilization using periodic terminal costs and sets with recursive feasibility and limit cycle tracking error asymptotic stability guarantees; 
\item[\textit{(ii)}] Systematic methods for computing quadratic periodic terminal costs and ellipsoidal or polytopic periodic terminal sets that satisfy the developed limit cycle stability and feasibility conditions for FCS-MPC.
\end{itemize}
Additional contributions include  formulation of an optimization problem for synthesis of an optimal limit cycle for switched affine systems with state constraints and a tractable procedure for computation of an outer polytopic approximation of the set of feasible states for limit cycle FCS-MPC. The set-iteration algorithm for computing \textit{polytopic} periodic invariant sets is a contribution beyond FCS-MPC, as it provides a larger admissible domain of attraction for switched affine systems compared to ellipsoidal invariant sets (in the case of polytopic constraints).

\begin{rem}
The FCS-MPC problem considered in this paper differs from the model predictive control problems for switched systems with dwell time constraints in previous studies, see for example \cite{zhuang2023model,chen2022efficient} and the references therein. The switched systems analyzed in these frameworks incorporate an additional continuous control set control input, which, together with dwell time constraints, enable stabilization of the switched system with respect to a fixed steady state point.
\end{rem}

\begin{rem}
In the continuous control set MPC (CCS-MPC) framework, periodic reference tracking problems have been studied in several papers \cite{limon2015mpc, aydiner2016periodic, kohler2018mpc, kohler2019nonlinear}. These papers cover both linear and nonlinear system classes and offer an extensive range of solutions for dealing with unreachable references and computing stabilizing terminal ingredients for tracking periodic references. 
Although the limit cycle represents a periodic reference, the CCS-MPC solutions for tracking periodic references use terminal control laws that take values in a continuous control set and thus, they cannot be used to establish asymptotic stability in the presence of a finite control set. 
\end{rem}

The remainder of this paper is organized as follows. 
Section~\ref{Section2} introduces the considered class of switched affine systems, conditions for existence of a limit cycle and methods for computing an optimal limit cycle. Section~\ref{Section3} describes the proposed FCS-MPC design and the limit cycle stability and feasible guarantees. Section~\ref{Section4} presents the algorithms to compute the terminal ingredients based on both ellipsoidal and polytopic set representations. Section~\ref{Section5} illustrates the implementation and effectiveness of limit cycle FCS-MPC for two different benchmarks. The conclusions are presented in Section~\ref{Section6}.

\paragraph*{Notation and basic definitions} The sets of real, non-negative real, integer and non-negative integer numbers are denoted as $\mathbb{R}$, $\mathbb{R}_{+}$, $\mathbb{I}$ and $\mathbb{I}_{+}$ respectively. $\mathbb{I}_{[i,j]}$ denotes the set of integer numbers restricted in $[i,j]$ ($i < j$).
The modulo operator is defined as $c$ = $a$ mod $b$ where $c$ is the remainder of the Euclidean division between the integers $a$ and $b$.
$\| x \|^2_P$ denotes the quadratic form $x^\top Px$.
The positive definite, positive semi-definite and negative semi-definite matrices are denoted by $A\succ0$, $A\succeq0$ and $A\preceq0$ respectively.
The notations $\det(A)$, $\lambda_{\min}(A)$ and $\lambda_{\max}(A)$ represent the determinant, the minimum and maximum eigenvalues of a matrix $A$.
$H_{i\bullet}$ denotes the $i$-th row of a matrix $H$.
Given two sets $\mathbb{X} \subseteq \mathbb{R}^n$ and $\mathbb{Y} \subseteq \mathbb{R}^n$,  the Minkowski set addition is defined by $\mathbb{X}\oplus\mathbb{Y}:=\{x+y\;|\;x\in\mathbb{X},y\in\mathbb{Y}\}$ 
and the Pontryagin set difference is defined by $\mathbb{X}\ominus\mathbb{Y}:=\{x\in\mathbb{X}\;|\;x+y\in\mathbb{X},\forall y\in\mathbb{Y}\}.$
$\text{int}(\mathbb{X})$ denotes the interior of the set $\mathbb{X}$. 
A real-valued scalar function $\phi:\mathbb{R}_{+}\rightarrow\mathbb{R}_{+}$ belongs to class $\mathcal{K}\; (\phi\in\mathcal{K})$ if it is continuous, strictly increasing and $\phi(0)=0$. 
A real-valued scalar function $\phi:\mathbb{R}_{+}\rightarrow\mathbb{R}_{+}$ belongs to class $\mathcal{K}_\infty\; (\phi\in\mathcal{K}_\infty)$ if $\phi\in\mathcal{K}$ and it is radially unbounded (i.e. $\phi(s)\rightarrow\infty$ as $s\rightarrow\infty$).

\section{Preliminaries} \label{Section2}

This paper considers the following class of discrete-time switched affine systems or, equivalently, the class of billinear affine systems with a finite control set, i.e.,
\begin{subequations} \label{equ2_switched_model_dis}
\begin{align}
x(k+1)&= A(u(k)) x(k) + b(u(k)), \label{equ2_switched_model_dis_1}\\
y(k)  &= C(u(k)) x(k) + d(u(k)),
\end{align}
\end{subequations}
where
$x \in \mathbb{R}^{n_x}$ is the system state, $u \in \mathbb{U}\subset \mathbb{R}^{n_u}$ is the control input and $y \in \mathbb{R}^{n_y}$ is the system output. The state and input are subject to the following constraints:
\begin{subequations} \label{}
\begin{align}
& x \in \mathbb{X} \subseteq \mathbb{R}^{n_x},\quad u \in \mathbb{U} \subset \mathbb{R}^{n_u},\quad y\in\mathbb{R}^{n_y} \\
&\hspace{-0.52em}\mathbb{U} := \left\{ u_1,\ldots,u_{N_s} \;|\; u_i \in \mathbb{R}^{n_u},\; \forall i \in \mathbb{I}_{[1,N_s]} \right\},
\end{align}
\end{subequations}
where $\mathbb{X}$ is a convex polytopic set, $\mathbb{U}$ is a \textit{finite control set} and $N_s$ represents the total number of elements of $\mathbb{U}$. For example, $\mathbb{U}=\{0,1\}$ is an often encountered case in power electronics to denote the status of the transistor (or of a binary switch) or, more generally, $\mathbb{U}=\{-0.2,-0.1,0,0.1,0.2\}$ can denote a quantized input signal within $[-0.2,0.2]$. Besides state constraints, also output constraints can be included, but we opted to omit output constraints to simplify the exposition.

One important real-life application of switched systems that can be modeled using \eqref{equ2_switched_model_dis} are switched-mode power converters, which are composed of electronic passive components. These components are operated by fast switching among a finite number of switch configurations \cite{spinu2012observer}:
\begin{subequations} \label{equ2_switched_model_con}
\begin{align}
\dot{x}(t)&= A_c(u(t)) x(t) +B_c(u(t))\omega, \\
y(t) &= C_c(u(t)) x(t) + D_c(u(t))\omega.
\end{align}
\end{subequations}
In this paper, a constant exogenous input $\omega$ is considered, which is typically corresponding to constant voltage and current sources for DC power converter applications. Thus $B_c(u(t))\omega$ and $D_c(u(t))\omega$ can be simplified as $b_c(u(t))$ and $d_c(u(t))$, respectively. The discrete-time switched affine system \eqref{equ2_switched_model_dis} can be derived by discretizing each corresponding subsystem \eqref{equ2_switched_model_con} for each switch mode using the zero-order-hold method for a suitable sampling time $T_s\in\mathbb{R}$.

\subsection{Limit cycle existence and computation}
The considered discrete-time switched affine system \eqref{equ2_switched_model_dis} admits a time-varying steady-state performance, which tends to have a natural convergence to a repetitive sequence, behaving as a limit cycle in discrete-time.
The limit cycle represents the stationary state of sustained oscillations, that depend exclusively on the parameters of the system and their intrinsic properties \cite{serieye2023attractors}.

To formally define a limit cycle for system \eqref{equ2_switched_model_dis} consider an arbitrary input sequence with length $p$, i.e., $U_{lc}=\{u_{lc}(0) ,\ldots , u_{lc}(p-1)\}$ and $u_{lc}(j)\in\mathbb{U}$. For brevity, the following notation is introduced:
\begin{align}
A_{j} := A(u_{lc}(j)), \quad
b_{j} := b(u_{lc}(j)).
\end{align}
A $p$-periodic state limit cycle $X_{lc}:=\{x_{lc}(0) ,\ldots , x_{lc}(p-1)\}$ with $x_{lc}(p)=x_{lc}(0)$ can then be derived by solving the following system of equations:
\begin{equation} \label{}
\hspace{-0.5em}
\begin{cases}
x_{lc}(1) = A_0 x_{lc}(0) + b_0, \\
x_{lc}(2) = A_1 x_{lc}(1) + b_1, \\
\hspace{3em}  \vdots \\
x_{lc}(0) = A_{p-1} x_{lc}(p-1) + b_{p-1}.
\end{cases}
\end{equation}
By defining $X_{p}:=(x_{lc}(0)^\top,\ldots,x_{lc}(p-1)^\top)^\top$ and the matrices:
\[M_p:=\begin{pmatrix}
  A_{0} & -I&      0 & \ldots & 0 \\
      0 &  A_{1} &     -I & \ddots & \vdots\\
 \vdots & \vdots & \ddots & \ddots & 0 \\
      0 &      0 & \ldots & A_{p-2} & -I \\
-I &      0 & \ldots &      0 & A_{p-1}
\end{pmatrix},\, \mathbf{b}_p:=\begin{pmatrix}
-b_0 \\ -b_1 \\ \vdots \\ -b_{p-2} \\ -b_{p-1}
\end{pmatrix},\]
we can compactly write the limit cycle equations as:
\begin{equation} \label{equ3_limit_cycle_calculation}
M_pX_p=\mathbf{b}_p.
\end{equation}
As implied by \eqref{equ3_limit_cycle_calculation}, the steady-state limit cycle exists if $M_p$ is invertible and a limit cycle sequence $X_{lc}$ can be directly obtained from the limit cycle vector $X_p$. In \cite{serieye2023attractors}, a necessary and sufficient condition for the existence of the limit cycle is provided, which is recalled next.

\begin{lem}\textit{\cite{serieye2023attractors}}
\label{lem1}
A chosen input sequence $U_{lc}$ with length $p$ generates a unique steady-state limit cycle for system \eqref{equ2_switched_model_dis} if and only if 1 is not an eigenvalue of the monodromy matrix $\mathcal{A}_p:=A_0\ldots A_{p-1}$. Moreover, the state limit cycle vector is given by $X_{p}=M_p^{-1}\mathbf{b}_p$ and the corresponding limit cycle sequence is $X_{lc}=\{x_{lc}(0) ,\ldots , x_{lc}(p-1)\}$. 
\end{lem}

As demonstrated for the system \eqref{equ2_switched_model_dis}, a $p$-periodic input sequence $U_{lc}$ can generate a unique corresponding $p$-periodic state limit cycle $X_{lc}$, which also characterizes the steady-state performance of the output $y(k)$. I.e., it yields a corresponding output trajectory $Y_{lc}:=\{y_{lc}(0) ,\ldots , y_{lc}(p-1)\}$, where $y_{lc}(j)$ is obtained from $x_{lc}(j)$ and $u_{lc}(j)$ via the system dynamics  \eqref{equ2_switched_model_dis}. 
In order to search for an optimal limit cycle with respect to a desired output reference and take the state constraints $\mathbb{X}$ into account, the following optimization problem can be formulated: 
\begin{subequations} \label{equ3_limit_cycle_problem}
\begin{align}
\min_{X_{lc},U_{lc}} &  W\left(X_{lc},U_{lc},\overline{Y}_p \right) \\
\text { s.t. } 
& x_{lc}(j+1)=A_j x_{lc}(j)+b_j, &&\hspace{-3em}\forall j \in \mathbb{I}_{[0,p-2]}, \label{equ3_limit_cycle_problem_ec1}\\
& x_{lc}(0) = A_{p-1}x_{lc}(p-1) + b_{p-1}, \label{equ3_limit_cycle_problem_ec2}\\
& u_{lc}(j) \in \mathbb{U},  &&\hspace{-3em}\forall j \in \mathbb{I}_{[0,p-1]},\\
& x_{lc}(j) \in \mathbb{X}, &&\hspace{-3em}\forall j \in \mathbb{I}_{[0,p-1]},
\end{align}
\end{subequations}
where $\overline{Y}_p=\{ \overline{y}(0), \ldots, \overline{y}(p-1) \}$ denotes a desired output reference of length $p$, which could be a constant or $p$-periodic reference. The cost function $W$ is defined as
\begin{equation} \label{equ3_limit_cycle_cost}
\begin{aligned}
W&(X_{lc},U_{lc},\overline{Y}_p):= \left\|\frac{1}{p}\sum_{j=0}^{p-1}\Big(y_{lc}(j)-\overline{y}(j)\Big)\right\|_1\\
&= 
\left\| \frac{1}{p} \sum_{j=0}^{p-1} \Big( C(u_{lc}(j)) x_{lc}(j)+d(u_{lc}(j))-\overline{y}(j) \Big)\right\|_1,
\end{aligned}
\end{equation}
which represents the $1$-norm of the mean error against the output reference at steady-state.
Alternatively, the cost function $W$ in \eqref{equ3_limit_cycle_cost} can be formulated using the 2-norm or $\infty$-norm of the mean error, to optimize a different steady-state output behavior \cite{geyer2016model,egidio2020global}. Note that if there are no state constraints, problem~\eqref{equ3_limit_cycle_problem} is always feasible under the hypothesis of Lemma~\ref{lem1}. In this paper we assume that this problem is feasible for the considered state constraints. 

\begin{assum} \label{assum_limit_cycle}
System \eqref{equ2_switched_model_dis} and the state constraints $\mathbb{X}$ are such that problem~\eqref{equ3_limit_cycle_problem} is feasible.
\end{assum}

For the remainder of the paper we use the notation $\overline{X}_{lc}:=\{\overline{x}_{lc}(0) ,\ldots , \overline{x}_{lc}(p-1)\}$ and 
$\overline{U}_{lc}:=\{\overline{u}_{lc}(0) ,\ldots , \overline{u}_{lc}(p-1)\}$ to denote the optimal steady-state and input limit cycle sequences calculated by solving problem~\eqref{equ3_limit_cycle_problem}. Note that problem~\eqref{equ3_limit_cycle_problem} is in general a mixed-integer nonlinear (bilinear) optimization problem. For many power electronics converter circuits the matrix $A$ does not depend on the switching input and then, problem~\eqref{equ3_limit_cycle_problem} is a mixed-integer linear program.

\section{Limit cycle FCS-MPC} \label{Section3}
To formulate the limit cycle FCS-MPC problem, the reference limit cycle corresponding to the optimal limit cycle sequence $\overline{X}_{lc}$ can be dynamically represented as a periodic solution of the system:
\begin{equation} \label{equ2_system_limit_cycle}
\begin{split}
\overline{x}(k+1) = \overline{A}_{k} \overline{x}(k) + \overline{b}_{k},
\end{split}
\end{equation}
where 
\begin{subequations} \label{}
\begin{align}
\overline{A}_{k} &:= A(\overline{u}(k)), \quad
\overline{b}_{k} := b(\overline{u}(k)), \\
\overline{x}(k)  &:= \overline{x}_{lc}(k \;\text{mod}\; p),\quad
\overline{u}(k)  := \overline{u}_{lc}(k \;\text{mod}\; p). 
\end{align}
\end{subequations}
Next, let $x_{0|k}=x(k)$ and $u_{0|k}=u(k)$ denote the system state and input at time $k$. Besides, let $x_{i|k}$ and $u_{i|k}$ denote the predicted state and input at time $k+i$ over a prediction horizon $N$. The sequences of the predicted state $X_k$ and the predicted inputs $U_k$ are defined as
\begin{subequations} \label{equ3_predicted_sequences}
\begin{align}
X_{k}&=\left\{ x_{0|k}, \ldots, x_{N|k} \right\},\\
U_{k}&=\left\{ u_{0|k}, \ldots, u_{N-1|k} \right\}.
\end{align}
\end{subequations}
If Assumption~\ref{assum_limit_cycle} is satisfied, then given the optimal state and input limit cycles $\overline{X}_{lc}$ and $\overline{U}_{lc}$, the time-varying cost function $J( x(k), U_{k}, k )$ of limit cycle FCS-MPC is defined as
\begin{equation} \label{equ3_FCSMPC_limit_cycle_J}
\begin{split}
J\left( x(k), U_{k}, k \right) :=& \sum_{i=0}^{N-1} l(x_{i|k}-\overline{x}_{i|k},u_{i|k}-\overline{u}_{i|k}) \\
&+ V_f(x_{N|k}-\overline{x}_{N|k},k),
\end{split}
\end{equation}
where $l:\mathbb{R}^{n_x}\times\mathbb{R}^{n_u}\rightarrow\mathbb{R}_{+}$ is the stage cost and $V_f:\mathbb{R}^{n_x}\times\mathbb{I}_{+}\rightarrow\mathbb{R}_{+}$ is a time-varying terminal cost.
$\overline{X}_{k}=\{\overline{x}_{0|k}, \ldots, \overline{x}_{N|k} \}$ and $\overline{U}_{k}=\{\overline{u}_{0|k}, \ldots, \overline{u}_{N-1|k} \}$ represent the sequences of the limit cycle dependent references, i.e.,
\begin{subequations} \label{equ4_mod}
\begin{align}
\overline{x}_{i|k} &= \overline{x}(k+i),  &&\forall i \in \mathbb{I}_{[0,N]},\\
\overline{u}_{i|k} &= \overline{u}(k+i),  &&\forall i \in \mathbb{I}_{[0,N-1]}.
\end{align}
\end{subequations}

Based on the cost function \eqref{equ3_FCSMPC_limit_cycle_J}, the system model \eqref{equ2_switched_model_dis} and the system constraints ($x(k) \in \mathbb{X}$, $u(k) \in \mathbb{U}$), the limit cycle FCS-MPC algorithm can be formulated as the following mixed-integer optimization problem, assuming that the full state is measurable (or observable):
\begin{subequations} \label{equ3_FCSMPC_limit_cycle_problem}
\begin{align}
\hspace{-0.9em}\min_{U_{k}}\;\; & J\left( x(k), U_{k}, k \right) \\
\hspace{-0.9em}\text { s.t. } 
& x_{i+1|k}=A(u_{i|k})x_{i|k}+b(u_{i|k}),  &&\hspace{-0.5em}\forall i \in \mathbb{I}_{[0,N-1]},\\
& \overline{x}_{i+1|k}=A(\overline{u}_{i|k})\overline{x}_{i|k}+b(\overline{u}_{i|k}),  &&\hspace{-0.5em}\forall i \in \mathbb{I}_{[0,N-1]},\\
& x_{i|k} \in \mathbb{X},  &&\hspace{-0.5em}\forall i \in \mathbb{I}_{[1,N-1]}, \\
& u_{i|k} \in \mathbb{U},  &&\hspace{-0.5em}\forall i \in \mathbb{I}_{[0,N-1]}, \\
& x_{N|k} \in \mathbb{X}_T(k), \label{equ3_terminal_condition}
\end{align}
\end{subequations}
where $\mathbb{X}_T(k)\subseteq\mathbb{R}^{n_x}$ is a time-varying terminal set.
For brevity, we denote
\begin{subequations} \label{equ3_notation_Ab}
\begin{align}
\overline{A}_{i|k}:=A(\overline{u}_{i|k}), \quad
\overline{b}_{i|k}:=b(\overline{u}_{i|k}).
\end{align}
\end{subequations}
The first optimal control action $u^\ast_{0|k}$ taking values in the finite control set $\mathbb{U}$ is used to define the FCS-MPC control law, i.e., $u(k) = \kappa_\text{mpc}(x(k)):=u^\ast_{0|k}$.

The problem considered in this paper is how to design the terminal cost $V_f$ and set $\mathbb{X}_T$ to guarantee recursive feasibility and asymptotic stability with respect to the limit cycle $\overline{X}_{lc}$ in the presence of a \textit{finite control set}.

\subsection{Terminal ingredients design}
Consider a set of $p$-periodic terminal ingredients $\{\kappa_j(\cdot), \mathcal{X}_{j} , F_j(\cdot)  \}$ with $j \in \mathbb{I}_{[0,p-1]}$, where $\kappa_j:\mathbb{R}^{n_x} \rightarrow \mathbb{U}$ denote the time dependent terminal control laws, $\mathcal{X}_{j}$ represent the terminal sets and $F_j:\mathbb{R}^{n_x} \rightarrow \mathbb{R}_+$ are the terminal costs.

The terminal $p$-periodic control laws are chosen as 
\begin{equation} \label{equ3_terminal_law}
\kappa_j(x(j)): = \overline{u}(j)\in\mathbb{U},\quad j \in \mathbb{I}_{[0,p-1]},
\end{equation}
i.e., they are taken equal to the elements of the optimal steady-state input limit cycle sequence $\overline{U}_{lc}$. For brevity, we use $\kappa_j$ to denote $\kappa_j(x(j))$. Next define $\Phi(x(j)):=A(\kappa_j) x(j) + b(\kappa_j)$.
When the switched affine system \eqref{equ2_switched_model_dis} is regulated by the terminal control law $\kappa_j$, the following switched affine autonomous system can be obtained:
\begin{equation} \label{equ3_system_xi}
\begin{aligned}
x(j+1) = \Phi(x(j)) = \overline{A}_j x(j) + \overline{b}_j.
\end{aligned}
\end{equation}

\begin{defn}[$p$-periodic invariant tube] \label{definition_invariant_tube} 
A sequence of $p$ sets $\mathcal{X}^\text{inv} = \{ \mathcal{X}_0, \ldots, \mathcal{X}_{p-1} \;:\; \mathcal{X}_i \subseteq \mathbb{X}, \; \forall j \in \mathbb{I}_{[0,p-1]}\}$ such that $\overline{x}(j) \in \text{int}(\mathcal{X}_j)$ is a $p$-periodic invariant tube for the switched affine autonomous system \eqref{equ3_system_xi} if for all 
$x(j) \in \mathcal{X}_{j},\; j \in \mathbb{I}_{[0,p-1]}$, it holds that $\Phi(x(j)) \in \mathcal{X}_{j+1 \;\text{mod}\; p}$.
\end{defn}

Note that the $p$-periodic invariant tube conditions can be alternatively formulated as
\begin{equation} \label{equ3_terminal_set}
\overline{A}_j \mathcal{X}_{j} \oplus \overline{b}_j \subseteq \mathcal{X}_{j+1 \;\text{mod}\; p},\; \forall j \in \mathbb{I}_{[0,p-1]}. 
\end{equation}
Next, consider the following conditions for designing a set of $p$-periodic terminal costs:
\begin{equation} \label{equ3_terminal_cost}
\begin{split}
F_{j+1 \;\text{mod}\; p}\left( \Phi_j(x(j)) -\overline{x}(j+1) \right) 
&- F_j\left( x(j)-\overline{x}(j) \right) \\
&\leq - l\left( x(j)-\overline{x}(j),0 \right), 
\end{split}
\end{equation}
for all $x(j) \in \mathcal{X}_{j}$ and $j \in \mathbb{I}_{[0,p-1]}$. Given $p$-periodic terminal sets and costs that satisfy the above conditions, the time-varying terminal set in the FCS-MPC problem \eqref{equ3_limit_cycle_problem} can be defined as 
\begin{equation} \label{equ3_MPC_terminal_set}
\begin{split}
\mathbb{X}_T(k) := \mathcal{X}_{k+N \;\text{mod}\; p},
\end{split}
\end{equation}
and the corresponding time-varying terminal cost can be defined as
\begin{equation} \label{equ3_MPC_terminal_cost}
\begin{split}
V_f(x_{N|k}-\overline{x}_{N|k},k) := F_{k+N \;\text{mod}\; p}\left( x_{N|k}-\overline{x}_{N|k} \right).
\end{split}
\end{equation}
Let $\mathbb{X}_f(N)$ denotes the set of feasible states, i.e., the set of all initial states that can admissibly reach the invariant tube $\mathcal{X}^\text{inv}$ in $N$ steps.
\begin{rem}
The proposed FCS-MPC problem \eqref{equ3_limit_cycle_problem} admits independent value selections concerning the limit cycle periodicity $p$ and the prediction horizon $N$, therefore, the time-varying terminal cost and set inherit a periodicity, which is the least common multiple of $p$ and $N$.
\end{rem}

\subsection{Recursive feasibility and limit cycle stability}
In this section, we study the closed-loop properties of the proposed limit cycle FCS-MPC algorithm. In particular, FCS-MPC recursive feasibility and closed-loop asymptotic stability of the limit cycle are analyzed. 

Stability of the discrete-time switched affine systems \eqref{equ2_switched_model_dis} with respect to a $p$-periodic time-varying state limit cycle can be analyzed by defining corresponding error dynamics and analyzing the asymptotic stability of the error with respect to the origin. When an analytic control law is available, the error dynamics can be defined explicitly as in \cite{serieye2023attractors}. In the case of MPC closed-loop systems, since the MPC control is computed numerically, the error dynamics can only be defined implicitly, as in  \cite{limon2015mpc}. Hence, in what follows we adopt the approach of \cite{limon2015mpc} by defining the tracking error $z(k):=x(k)-\overline{x}(k)$ and its implicit autonomous dynamics 
\begin{equation} \label{equ3_error_dynamic}
z(k+1) =f_\text{mpc}(z(k)),\quad k\in\mathbb{I}_+.
\end{equation}

\begin{assum} \label{assum_costs} 
The stage cost $l(\cdot,\cdot)$ and the terminal cost $V_f(\cdot,\cdot)$ satisfy the following conditions:
\begin{subequations} \label{equ3_cost_bounds}
\begin{align}
l(x(k)-\overline{x}(k),u(k)-\overline{u}(k)) &\geq \alpha_1(\|x(k)-\overline{x}(k)\|), \nonumber\\
&\forall x(k) \in \mathbb{X}_f(N), \; \forall u(k) \in \mathbb{U}, \label{equ3_cost_bounds_a}\\
V_f(x(k)-\overline{x}(k),k) &\leq \beta_1(\|x(k)-\overline{x}(k)\|), \nonumber\\
&\forall x(k) \in \mathcal{X}_{k\;\text{mod}\;p},\label{equ3_cost_bounds_b}
\end{align}
\end{subequations}
where $\alpha_1(\cdot)$ and $\beta_1(\cdot)$ are $\mathcal{K}_\infty$ functions.
\end{assum}
In what follows we state a result that combines \cite[Theorem~1]{limon2015mpc} and \cite[Theorem~2]{jiang2002converse}, as limit cycle FCS-MPC uses a time-varying cost function.
\begin{prop}\cite{limon2015mpc, jiang2002converse} \label{theorem_Lyapunov}  
Consider the FCS-MPC closed-loop error dynamics \eqref{equ3_error_dynamic}. Let $\Gamma\subseteq \mathbb{R}^{n_x}$ be a positively invariant set for the error dynamics \eqref{equ3_error_dynamic}, which contains the origin in its interior. Suppose that there exists a function $V:\mathbb{R}^{n_x}\times \mathbb{I}_+ \rightarrow \mathbb{R}_{+}$ and suitable $\mathcal{K}_\infty$-class functions $\alpha_1$, $\alpha_2$, $\alpha_3$ such that
\begin{equation}
\begin{split}
\textit{(i)}~&V(z(k),k)\geq \alpha_1(\| z(k) \|),\; \forall z(k)\in \Gamma, \\
\textit{(ii)}~&V(z(k),k)\leq \alpha_2(\| z(k) \|),\; \forall z(k)\in \Gamma, \\
\textit{(iii)}~&V(z(k+1),k+1)-V(z(k),k) \\
& \hspace{4.5em} \leq - \alpha_3(\| z(k) \|),\; \forall z(k)\in \Gamma.
\end{split}
\end{equation}
Then $V(\cdot,\cdot)$ is called a (time-varying) Lyapunov function and the origin of the error dynamics \eqref{equ3_error_dynamic} is asymptotically stable for all $z(0)$ in $\Gamma$.
\end{prop}

Next consider the limit cycle tracking FCS-MPC problem \eqref{equ3_FCSMPC_limit_cycle_problem}.  
Let $U_{k}^* = \{u_{0|k}^{*}, \ldots, u_{N-1|k}^{*}\}$ denote the optimal input sequence at time $k$ and 
let $X_{k}^* = \{x_{0|k}^{*}, \ldots, x_{N|k}^{*}\}$ denote optimal state sequence.

\begin{thm} \label{theorem_RF}
Consider the limit cycle FCS-MPC problem \eqref{equ3_FCSMPC_limit_cycle_problem}. 
Let the time-varying terminal set be defined as in \eqref{equ3_MPC_terminal_set} and assume that $\mathcal{X}^\text{inv}$ is a $p$-periodic invariant tube as in Definition~\ref{definition_invariant_tube}. 
Let the time-varying terminal cost be defined as in \eqref{equ3_MPC_terminal_cost} with $\{F_j\}_{j\in\mathbb{I}_{[0,p-1]}}$ satisfying the conditions \eqref{equ3_terminal_cost}.
Let the stage cost and terminal cost satisfy Assumption~\ref{assum_costs}.
Then, the following statements hold:
\begin{enumerate}
\item[\textit{(i)}]~For all $x(0)\in \mathbb{X}_f(N)$, the FCS-MPC problem \eqref{equ3_FCSMPC_limit_cycle_problem} is feasible for all $k\in \mathbb{I}_+$. \\
\item[\textit{(ii)}]~The origin of the closed-loop FCS-MPC error dynamics \eqref{equ3_error_dynamic} is asymptotically stable and $x(k)$ converges asymptotically to $\overline{x}(k)$ for all $x(0)\in \mathbb{X}_f(N)$.
\end{enumerate}
\end{thm}

\begin{pf}
\textit{(i)}~Recursive feasibility proof: Since $x(0)\in \mathbb{X}_f(N)$, the FCS-MPC problem \eqref{equ3_FCSMPC_limit_cycle_problem} is feasible at time $k=0$. Then we proceed by induction, i.e., suppose that the FCS-MPC problem \eqref{equ3_FCSMPC_limit_cycle_problem} is feasible at time $k$ and define a shifted, sub-optimal input sequence at time $k+1$:
\begin{equation} \label{equ3_shifted_sequence}
\begin{split}
\widetilde{U}_{k+1} &= \left\{\widetilde{u}_{0|k+1}, \ldots, \widetilde{u}_{N-2|k+1},\widetilde{u}_{N-1|k+1}\right\} \\
&:=\left\{u_{1|k}^{*}, \ldots, u_{N-1|k}^{*},\widetilde{u}_{N-1|k+1}\right\}.
\end{split}
\end{equation}
Constraint \eqref{equ3_terminal_condition} implies that $x_{N|k}^{*}\in \mathcal{X}_{k+N \;\text{mod}\; p}$; thus the terminal control law \eqref{equ3_terminal_law} can be used as $\widetilde{u}_{N-1|k+1}$, i.e.,
\begin{equation}
\begin{split}
\widetilde{u}_{N-1|k+1} := \kappa_{N+k \;\text{mod}\; p}(x_{N|k}^{*})= \overline{u}_{N|k} \in \mathbb{U}.
\end{split}
\end{equation}
Following the notation \eqref{equ3_notation_Ab}, we can infer that:
\begin{subequations} \label{equ3_shifted_law}
\begin{align}
\widetilde{x}_{i-1|k+1} &= x_{i|k}^{*} \in \mathbb{X}, \; \forall i \in \mathbb{I}_{[1,N]}, \\
\widetilde{x}_{N|k+1} &= \overline{A}_{N|k}x_{N|k}^*+\overline{b}_{N|k} \in \mathcal{X}_{k+N+1 \;\text{mod}\; p}.
\end{align}
\end{subequations}
According to \eqref{equ3_MPC_terminal_set}, it can be inferred that $\widetilde{x}_{N|k+1} \in \mathbb{X}_T(k+1)$.
Hence the FCS-MPC problem \eqref{equ3_FCSMPC_limit_cycle_problem} remains feasible at time $k+1$, which completes the recursive feasibility proof.

\textit{(ii)}~Asymptotic stability proof: From the proof of recursive feasibility we have that $x(k)\in\mathbb{X}_f(N)$ for all $k\geq 0$. Thus, since the tracking error is defined as $z(k)=x(k)-\overline{x}(k)$, using the limit cycle dynamics \eqref{equ2_system_limit_cycle} we can define $\Gamma:=\cup_{j\in\mathbb{I}_{[0,p-1]}}\Gamma_j$ with $\Gamma_j:=\mathbb{X}_f(N)\ominus \overline{x}(j)$ such that $z(k)\in\Gamma$ for all $k\geq 0$ and $0\in\text{int}(\Gamma)$. Hence, $\Gamma$ is a positively invariant set for the FCS-MPC error dynamics \eqref{equ3_error_dynamic}. 

Next we will show that the optimal value function of the FCS-MPC problem
\begin{equation} \label{}
\begin{split}
V(z(k),k) &= V( x(k)-\overline{x}(k) , k )=J\left( x(k), U_{k}^*, k \right)
\end{split}
\end{equation}
is a time-varying Lyapunov function for the FCS-MPC error dynamics, i.e., it  satisfies the conditions of Proposition~\ref{theorem_Lyapunov} in $\Gamma$. For any $k\geq 0$, it holds that $z(k)\in\Gamma_{k\;\text{mod}\;p}$ and $x(k)=z(k)+\overline{x}(k)\in\mathbb{X}_f(N)$. 

From the definition of the optimal value function and the assumption on the stage cost \eqref{equ3_cost_bounds_a}, we have that
\begin{equation*} 
\begin{split}
V(z(k),k) \geq  l(x_{0|k}-\overline{x}_{0|k},&u_{0|k}-\overline{u}_{0|k}) \geq  \alpha_1(\| z(k) \|), \\
 &\,\forall z(k)\in\Gamma_{k\;\text{mod}\;p}\subset\Gamma,
\end{split}
\end{equation*}
which is condition \textit{(i)} in Proposition~\ref{theorem_Lyapunov}. The terminal cost satisfies the following upper bound condition \eqref{equ3_cost_bounds_b}:
\begin{equation} \label{}
\begin{split}
V_f(z(k),k) \leq & \beta_1( \| x(k)-\overline{x}(k) \| ), \,\forall x(k)  \in \mathcal{X}_{k\;\text{mod}\;p}.
\end{split}
\end{equation}
Hence, for any $\widetilde{x}_{0|k}=x(k) \in \mathcal{X}_{k\;\text{mod}\;p}$, the terminal control law \eqref{equ3_terminal_law} can be applied to generate a sub-optimal state sequence that satisfies (due to \eqref{equ3_MPC_terminal_cost}):
\begin{equation} \label{}
\begin{split}
V_f(\widetilde{x}_{i|k}-\overline{x}_{i|k},k+i) = F_{k+i\;\text{mod}\;p}(\widetilde{x}_{i|k}-\overline{x}_{i|k}), \\
\widetilde{x}_{i|k} \in \mathcal{X}_{k+i\;\text{mod}\;p}.
\end{split}
\end{equation}
Then, the terminal cost condition \eqref{equ3_terminal_cost} implies that
\begin{equation} \label{}
\begin{split}
V_f(\widetilde{x}_{1|k}-\overline{x}_{1|k},k+1) + l(\widetilde{x}_{0|k}-\overline{x}_{0|k},0)&\\
\leq V_f(\widetilde{x}_{0|k}-\overline{x}_{0|k},k) &,\\
V_f(\widetilde{x}_{2|k}-\overline{x}_{2|k},k+2) + l(\widetilde{x}_{1|k}-\overline{x}_{1|k},0)& \\
\leq V_f(\widetilde{x}_{1|k}-\overline{x}_{1|k},k+1)&,\\
\ldots\\
V_f(\widetilde{x}_{N|k}-\overline{x}_{N|k},k+N) + l(\widetilde{x}_{N-1|k}-\overline{x}_{N-1|k},0) &\\
\leq V_f(\widetilde{x}_{N-1|k}-\overline{x}_{N-1|k},k+N-1)&.
\end{split}
\end{equation}
Summing up the above inequalities yields:
\begin{equation} \label{}
\begin{split}
V(&z(k),k) \leq J(x(k),\widetilde{U}_k,k) \\
&= \sum_{i=0}^{N-1} l(\widetilde{x}_{i|k}-\overline{x}_{i|k},0)+  V_f(\widetilde{x}_{N|k}-\overline{x}_{N|k},k+N) \\
&\leq V_f(\widetilde{x}_{0|k}-\overline{x}_{0|k},k),
\end{split}
\end{equation}
which implies
\begin{equation} \label{}
\begin{split}
V(z(k),k) \leq& \beta_1( \| z(k) \| ), \\
&\forall z(k) \in \mathcal{X}_{k\;\text{mod}\;p}\ominus\overline{x}(k) \subseteq \Gamma_{k\;\text{mod}\;p}.
\end{split}
\end{equation}
Next, we show that the upper bound condition can be extended to the positively invariant set $\Gamma$.
Since $\overline{x}(j) \in \text{int}(\mathcal{X}_j)$ as defined in Definition~\ref{definition_invariant_tube}, it gives that $0\in\text{int}(\Gamma_j)$ and we can find a $r>0$ such that
\begin{equation} \label{}
\begin{split}
\mathcal{B}_{r}:=\{z\in\mathbb{R}^{n_x} \;|\; \| z \|\leq r\} \subseteq \cap_{j\in\mathbb{I}_{[0,p-1]}}\{\mathcal{X}_{j}\ominus\overline{x}(j)\}.
\end{split}
\end{equation}
In addition, due to the boundedness of the state and input constraints, there exists a $d>0$ such that $V(z(k),k)\leq d$ for all $z(k)  \in \Gamma_{k\;\text{mod}\;p}$, which implies that for all $z(k)\in\Gamma_{k\;\text{mod}\;p}\setminus\mathcal{B}_{r}$, due to $\|z(k)\|>r$, 
\[V(z(k),k)\leq \beta_2(\|z(k)\|):=\frac{d}{r}\|z(k)\|.\]
Since $V(z(k),k) \leq \beta_1(\|z(k)\|)$ for all 
$z(k)\in\mathcal{B}_{r}$, we obtain that 
\begin{equation*}
\begin{split}
V(z(k),k) &\leq \alpha_2( \| z(k) \| ),\; \forall z(k)  \in \Gamma , \\
\alpha_2(\cdot) &= \max\{  \beta_1(\cdot) , \beta_2(\cdot) \}\in\mathcal{K}_\infty,
\end{split}
\end{equation*}
which is condition~\textit{(ii)} in Proposition~\ref{theorem_Lyapunov}.

Next, considering again the sub-optimal but feasible input and state sequences \eqref{equ3_shifted_sequence}-\eqref{equ3_shifted_law}, it holds that
\begin{equation} \label{equ3_Delta_V}
\begin{split}
& V(z(k+1),k+1) - V(z(k),k) \\
&= J( x(k+1), U_{k+1}^*, k+1 ) - J( x(k), U_{k}^*, k )\\
&\leq J( x(k+1), \widetilde{U}_{k+1}, k+1 ) - J( x(k), U_{k}^*, k )\\
&=  V_f(\widetilde{x}_{N|k+1}-\overline{x}_{N|k+1},k+1) \\
 & -V_f(x_{N|k}^*-\overline{x}_{N|k},k) + l(x_{N|k}^*-\overline{x}_{N|k},0) \\
 & - l(x_{0|k}-\overline{x}_{0|k},u_{0|k}-\overline{u}_{0|k}).
\end{split}
\end{equation}
Then according to \eqref{equ3_terminal_cost}, \eqref{equ3_MPC_terminal_cost} and \eqref{equ3_Delta_V} we obtain
\begin{equation} \label{equ3_Lyapunov3}
\begin{split}
 V(z(k+1),k+1) &- V(z(k),k) \\&\leq - l(x_{0|k}-\overline{x}_{0|k},u_{0|k}-\overline{u}_{0|k}) \\
&\leq - \alpha_1(\| z(k) \|), \; \forall z(k)  \in \Gamma,
\end{split}
\end{equation}
which is condition \textit{(iii)} in Proposition~\ref{theorem_Lyapunov}. Thus, the optimal value function $V(z(k),k) = V( x(k)-\overline{x}(k) , k )$ is a time-varying Lyapunov function for the FCS-MPC closed-loop dynamics and the origin is an asymptotically stable equilibrium of the limit cycle FCS-MPC error dynamics. In turn, this implies that the closed-loop state trajectory $x(k)$ asymptotically converges to the optimal limit cycle $\overline{x}(k)$.
\hfill $\qed$
\end{pf}

\section{Computation of the terminal ingredients} \label{Section4}
In this section we provide systematic methods for computing periodic terminal costs and sets for limit cycle FCS-MPC that satisfy the assumptions required for recursive feasibility and asymptotic stability. 
\subsection{Computation of the terminal costs}
In this paper, we consider quadratic cost functions, i.e.,
\begin{equation} \label{equ4_cost_l}
\begin{split}
l(x_{i|k}-\overline{x}_{i|k},u_{i|k}-\overline{u}_{i|k}) = &
\| x_{i|k}-\overline{x}_{i|k}\|^2_Q \\
& + \| u_{i|k}-\overline{u}_{i|k}\|^2_R,
\end{split}
\end{equation}
where $Q$ and $R$ are positive definite matrices.
Correspondingly, the time-varying terminal costs are formulated as:
\begin{equation} \label{equ4_cost_F}
\begin{split}
V_f(x_{N|k}-\overline{x}_{N|k},k) = \| x_{N|k}-\overline{x}_{N|k} \|^2_{P_{N|k}},
\end{split}
\end{equation}
where 
\begin{equation} \label{equ4_cost_P}
\begin{split}
P_{N|k} &= P_{lc}(k+N\; \text{mod}\; p)
\end{split}
\end{equation}
and $\mathbb{P}_{lc}=\{P_{lc}(0),\ldots,P_{lc}(p-1)\}$ denotes a set of $p$-periodic positive definite matrices.

Consider the quadratic cost functions \eqref{equ4_cost_l}-\eqref{equ4_cost_F}, and the following linear matrix inequalities (LMIs) with variables $Q$ and $\mathbb{P}_{lc}$:
\begin{subequations} \label{equ4_P_property}
\begin{align}
Q \succ 0, \; P_{lc}(j)\succ 0, \;
\forall j \in \mathbb{I}_{[0,p-1]}. \label{equ4_P_property_a} \\
\overline{A}_j^\top P_{lc}(j+1 \;\text{mod}\; p) \overline{A}_j -P_{lc}(j)+Q \preceq 0, \nonumber\\
\forall j \in \mathbb{I}_{[0,p-1]}. \label{equ4_P_property_b}
\end{align}
\end{subequations} 

\begin{lem} \label{lemma_terminal_costs}
Let the solutions $Q$ and $\mathbb{P}_{lc}$ of \eqref{equ4_P_property} be used in the cost functions \eqref{equ4_cost_l}-\eqref{equ4_cost_F}. Then the cost functions \eqref{equ4_cost_l}-\eqref{equ4_cost_F} satisfy Assumption~\ref{assum_costs} and the terminal costs satisfy condition \eqref{equ3_terminal_cost}.
\end{lem}

\begin{pf}
For the stage cost function $l(\cdot,\cdot)$ \eqref{equ4_cost_l} with positive definite matrices $Q$ and $R$ it holds that
\begin{equation}
\begin{split}
l(x(k)-\overline{x}(k),u(k)-\overline{u}(k)) & \geq \| x(k)-\overline{x}(k) \|^2_Q \\
& \geq \lambda_{\min}(Q) (\| x(k)-\overline{x}(k) \|) ,\\
&\forall x(k) \in \mathbb{X}_f(N), \; \forall u(k) \in \mathbb{U},
\end{split}
\end{equation}
where $\alpha_1(\cdot):=\lambda_{\min}(Q) (\cdot)$ is a $\mathcal{K}_\infty$ function and \eqref{equ3_cost_bounds_a} in Assumption~\ref{assum_costs} is satisfied.
Similarly, for the terminal cost function $V_f(\cdot,\cdot)$ \eqref{equ4_cost_F} defined as in \eqref{equ4_cost_P} it holds that
\begin{equation}
\begin{split}
V_f(x(k)-\overline{x}(k),k) & = \| x(k)-\overline{x}(k) \|^2_{P_{N|k}} \\
& \leq \lambda_{\max}(P_{N|k}) (\| x(k)-\overline{x}(k) \|)\\
&\leq \Lambda_{\max} (\| x(k)-\overline{x}(k) \|),\\
&\forall x(k) \in \mathcal{X}_{k\;\text{mod}\;p},
\end{split}
\end{equation}
where $\Lambda_{\max}:=\max\{ \lambda_{\max}(P_{lc}(j)),\;\forall j \in \mathbb{I}_{[0,p-1]}\}>0$ and \eqref{equ3_cost_bounds_b} in Assumption~\ref{assum_costs} holds.

Next, consider the quadratic cost functions \eqref{equ4_cost_l}-\eqref{equ4_cost_F}, i.e.,
\begin{subequations} \label{equ4_proof1}
\begin{align}
l(x(j)-\overline{x}(i),0) &= \| x(j) - \overline{x}(j) \|^2_Q, \\
F_j(x(j)-\overline{x}(i)) &= \| x(j) - \overline{x}(j) \|^2_{P_{lc}(j)},
\end{align}
\end{subequations} 
for all $x(j) \in \mathcal{X}_{j}$ and $j \in \mathbb{I}_{[0,p-1]}$. Recalling the autonomous systems \eqref{equ3_system_xi}$, i.e., \Phi_i(x(j)) = \overline{A}_j x(j) + \overline{b}_j$  and \eqref{equ2_system_limit_cycle}, i.e., $\overline{x}(j+1) = \overline{A}_j \overline{x}(j) + \overline{b}_j$, we obtain
\begin{equation} \label{equ4_proof2}
\begin{split}
F_{j+1 \;\text{mod}\; p} ( \Phi_j (x(j)) &-\overline{x}(j+1) ) \\
&= \| x(j) - \overline{x}(j) \|^2_{P_{lc}(j+1 \;\text{mod}\; p)}, \\
& \forall x(j) \in \mathcal{X}_{j} \; \forall j \in \mathbb{I}_{[0,p-1]}.
\end{split}
\end{equation} 
Thus, from \eqref{equ4_proof1}, \eqref{equ4_proof2} and\eqref{equ4_P_property_b} we obtain \eqref{equ3_terminal_cost}, which completes the proof.
\hfill $\qed$ 
\end{pf}

\subsection{Computation of the terminal sets}
Consider next the error dynamics corresponding to the closed-loop switched system \eqref{equ3_system_xi}, i.e.,
\begin{equation} \label{equ4_autonomous_error_dynamic}
	\begin{split}
		z(j+1)&= \overline{A}_j x(j) + \overline{b}_j
		- \overline{A}_j \overline{x}(j) - \overline{b}_j \\
		&=\overline{A}_j z(j), \; \forall z(j) \in \mathcal{X}_{j \;\text{mod}\; p}\ominus \overline{x}(j).
	\end{split}
\end{equation}
The next Lemma shows that invariant tubes for the error dynamics yield invariant tubes for the state dynamics. 
\begin{lem} \label{lem_Zset}
Let a feasible $p$-periodic invariant tube $\mathcal{Z}^\text{inv} = \{ \mathcal{Z}_0, \ldots, \mathcal{Z}_{p-1} \;|\; \mathcal{Z}_i \subseteq \mathbb{X} \ominus \overline{x}(j), \quad \forall j \in \mathbb{I}_{[0,p-1]}\}$ that contains the origin in each set $\mathcal{Z}_j$ be given for the autonomous error dynamics \eqref{equ4_autonomous_error_dynamic}.
Then the $p$-periodic tube $\mathcal{X}^\text{inv}$ constructed as
\begin{equation}
\mathcal{X}^\text{inv} = \{ \mathcal{X}_0, \ldots, \mathcal{X}_{p-1} \;:\; \mathcal{X}_{j} = \mathcal{Z}_{j} \oplus \overline{x}(j), \; \forall j \in \mathbb{I}_{[0,p-1]}\}
\end{equation}
is a $p$-periodic \textit{invariant} tube for the autonomous state dynamics \eqref{equ3_system_xi}.
\end{lem}
 
\begin{pf}
Considering that $z(j)\in \mathcal{Z}_{j} \subseteq \mathbb{X} \ominus \overline{x}(j)$, we obtain
\begin{equation}
\begin{split}
x(j) = z(j) + \overline{x}(j) \in \mathcal{Z}_{j} \oplus \overline{x}(j) := \mathcal{X}_{j} \subseteq \mathbb{X}.
\end{split}
\end{equation}
Similarly, due to \[z(j+1) = \overline{A}_j z(j) \in \mathcal{Z}_{j+1 \;\text{mod}\; p} \subseteq \mathbb{X} \ominus \overline{x}(j+1),\] it can be inferred that
\begin{equation}
\begin{split}
x(j+1) &= \overline{A}_j x(j) + \overline{b}(j) =  \overline{A}_j z(j) + \overline{A}_i \overline{x}(j) + \overline{b}(j) \\
&= \overline{A}_j z(j) + \overline{x}(j+1) \\
& \in \mathcal{Z}_{j+1 \;\text{mod}\; p} \oplus \overline{x}(j+1) := \mathcal{X}_{j+1 \;\text{mod}\; p}\subseteq \mathbb{X}.
\end{split}
\end{equation}
which completes the proof.
\hfill $\qed$
\end{pf}
Next, we present the procedure for computing an \textit{ellipsoidal $p$-periodic invariant tube}. To this end, consider $p$ ellipsoidal sets, which are defined as
\begin{equation} \label{equ4_ellip_sets}
\mathcal{Z}_{j} = \left\{  z \in \mathbb{R}^{n_x}  \; | \;  z^\top Z_j z\leq 1  \right\} , \quad \forall j \in \mathbb{I}_{[0,p-1]},
\end{equation}
where $Z_i$ are positive definite matrices and 
each set $\mathcal{Z}_{j}$ should satisfy the following constraints:
\begin{equation} \label{equ4_Z_constraint}
\mathcal{Z}_{j} \subseteq \mathbb{X} \ominus \overline{x}(j), \quad \forall j \in \mathbb{I}_{[0,p-1]}.
\end{equation}
Let the polytopic constraint set $\mathbb{Z}_{j} := \mathbb{X} \ominus \overline{x}(j)$ be explicitly parameterized as 
\begin{equation} \label{equ4_constraint_z_def}
\mathbb{Z}_{j} = \{  z \in \mathbb{R}^{n_x}  \; | \;  H_{\mathbb{Z}_{j}} \leq 1_{h_{\mathbb{Z}_{j}}} \}, \quad \forall j \in \mathbb{I}_{[0,p-1]},
\end{equation}
and consider the following tube-synthesis problem.
\begin{prob} \label{problem_ellipsoidal}
\begin{subequations} \label{equ4_LMI}
\begin{align}
\min_{O_j}\; & -\sum_{j=0}^{p-1} \log\det \left(O_j\right) \label{equ4_obj} \\
\text { s.t. } 
& O_j \succ 0 , \; \forall j \in \mathbb{I}_{[0,p-1]}, \label{equ4_con1} \\ 
& \begin{pmatrix}
O_j & O_j \overline{A}_j^\top \\ \overline{A}_j O_j & O_{j+1 \;\text{mod}\; p}
\end{pmatrix} \succeq 0,  \; j \in \mathbb{I}_{[0,p-1]}, \label{equ4_con2} \\
& \left\{ [H_{\mathbb{Z}_{j}}]_{i\bullet} O_j [H_{\mathbb{Z}_{j}}]_{i\bullet}^\top\leq 1, \forall i \in \mathbb{I}_{[1,h_{\mathbb{Z}_{j}}]} \right\}, 
\forall j \in \mathbb{I}_{[0,p-1]}. \label{equ4_con3}
\end{align}
\end{subequations}
\end{prob}

\begin{lem} \label{lemma_ellipsoidal}
Suppose there exists feasible solutions $\{O_j\}_{j \in \mathbb{I}_{[0,p-1]}}$ of Problem~\ref{problem_ellipsoidal} and let the $p$ ellipsoidal sets in \eqref{equ4_ellip_sets} be defined with $Z_j = O_j^{-1}$, for all $j \in \mathbb{I}_{[0,p-1]}$. Then the tube $\mathcal{Z}^\text{inv} = \{ \mathcal{Z}_0, \ldots, \mathcal{Z}_{p-1} \}$ with
\begin{equation}
\begin{split}
\mathcal{Z}_{j} = \{  z \in \mathbb{R}^{n_x}  \; | \;  z^\top Z_j z\leq 1  \}, \;
\forall j \in \mathbb{I}_{[0,p-1]},
\end{split}
\end{equation}
is a $p$-periodic invariant tube for the error dynamics \eqref{equ4_autonomous_error_dynamic}.
\end{lem}

\begin{pf}
Constraint \eqref{equ4_con1} implies that $Z_j$ is positive definite, i.e.,
\begin{equation} 
Z_j \succ 0 , \quad \forall j \in \mathbb{I}_{[0,p-1]}.
\end{equation}
Applying the Schur complement to \eqref{equ4_con2}, pre and post multiplying with $Z_j$ we obtain
\[Z_j-\overline{A}_j^\top Z_{j+1 \;\text{mod}\; p} \overline{A}_j\geq 0,  \; \forall j \in \mathbb{I}_{[0,p-1]},\]
which yields \cite{egidio2020global,serieye2023attractors}
\begin{equation}
z^\top \overline{A}_j^\top Z_{j+1 \;\text{mod}\; p} \overline{A}_j z\leq z^\top Z_j z, \; \forall j \in \mathbb{I}_{[0,p-1]}.\\
\end{equation}
From the above inequality we can directly deduce the invariance properties in Definition~\ref{definition_invariant_tube}. Finally, as shown in \cite{boyd1994linear}, \eqref{equ4_con3}  ensures that each ellipsoidal set is contained within the corresponding polytopic constraint set \eqref{equ4_Z_constraint}-\eqref{equ4_constraint_z_def}.
\hfill $\qed$
\end{pf}
Note that in order to maximize the volume of the ellipsoidal sets the cost \eqref{equ4_obj} is adopted as the objective function in the optimization problem \cite{boyd1994linear}.

Ellipsoidal tubes are of interest because the computation reduces to a convex optimization problem. However, such sets introduce convex quadratic terminal constraints in the FCS-MPC problem. Thus, next we propose a set-iterations algorithm for computing \textit{polytopic $p$-periodic invariant tubes}, which lead to linear terminal constraints. 

Let $\{ \mathcal{Z}_0, \ldots, \mathcal{Z}_{p-1}\}$ be a set of $p$ polytopes. Inspired by the standard set recursion algorithm for computing the invariant sets for linear systems \cite{kvasnica2004multi}, we propose a set recursion algorithm for computing a $p$-periodic invariant tube with polytopic sets, as summarized in Algorithm~\ref{algorithm_polytopic}. 

\begin{algorithm}[h]
\caption{Set recursion for computing a polytopic $p$-periodic invariant tube} \label{algorithm_polytopic}
\hspace*{\algorithmicindent} \textbf{Initialization:} $\mathcal{Z}_{j}^{(0)} \leftarrow \mathbb{X} \ominus \overline{x}(j), \;\forall j \in \mathbb{I}_{[0,p-1]}$
\begin{algorithmic}[1]
\For {$n=1:1:n_{\max}$}
    \For {$j=p-1:-1:0$}
    \If { $j = p$}
        \State $\mathcal{Z}_{p-1}^{(n)} \leftarrow \left\{ z \;|\; \overline{A}_{p-1} z \in \mathcal{Z}_{0}^{(n-1)} \right\}$;
    \Else
        \State $\mathcal{Z}_{j}^{(n)} \leftarrow \left\{ z \;|\; \overline{A}_j z \in \mathcal{Z}_{j+1}^{(n)} \right\}$;
    \EndIf
    \State $\mathcal{Z}_{j}^{(n)} \leftarrow \mathcal{Z}_{j}^{(n)} \bigcap \mathcal{Z}_{j}^{(n-1)}$;
    \EndFor
\EndFor
\end{algorithmic}
\end{algorithm}

The proposed algorithm operates the set recursion for a sufficient large number of iterations $n_{\max}$.
For each set $\mathcal{Z}_{j}^{(n)}$, $n$ denotes the iteration number and $\mathcal{Z}_{j}^{(0)} := \mathbb{X} \ominus \overline{x}(j)$, which is the corresponding polytopic constraint set as defined in \eqref{equ4_constraint_z_def}.
The set recursion is performed by computing the backward reachable sets for the autonomous error dynamics \eqref{equ4_autonomous_error_dynamic} as demonstrated in step 4 and 6.
Step 8 updates the derived backward reachable set $\mathcal{Z}_{j}^{(n)}$ by intersecting it with $\mathcal{Z}_{j}^{(n-1)}$ at the previous iteration, which ensures the constraints and promotes convergence of the iterations.
The proposed algorithm can be terminated earlier if $\mathcal{Z}_{j}^{(n)}=\mathcal{Z}_{j}^{(n-1)}$, for all $j \in \mathbb{I}_{[0,p-1]}$.

\begin{lem} \label{lemma_polytopic}
Suppose that Algorithm~\ref{algorithm_polytopic} converges in finite iterations to sets with non-empty interior and let $\mathcal{Z}_{j}^{(n_{\max})}$ be the resulting sets. Then the tube constructed as $\mathcal{Z}^\text{inv} = \{ \mathcal{Z}_0, \ldots, \mathcal{Z}_{p-1} \}$ with
\begin{equation}
\begin{split}
\mathcal{Z}_{j} = \mathcal{Z}_{j}^{(n_{\max})}, \;
\forall j \in \mathbb{I}_{[0,p-1]},
\end{split}
\end{equation}
is a $p$-periodic invariant tube for the error dynamics \eqref{equ4_autonomous_error_dynamic}.
\end{lem}

\begin{pf}
The finite iteration convergence assumption implies that
$\mathcal{Z}_{j}^{(n_{\max})} = \mathcal{Z}_{j}^{(n_{\max}-1)}$, for all $j \in \mathbb{I}_{[0,p-1]}$. As indicated by step 4 in Algorithm~\ref{algorithm_polytopic}, it follows:
\begin{equation}
\begin{split}
\mathcal{Z}_{p-1}^{(n_{\max})} 
&= \left\{ z \;|\; \overline{A}_{p-1} z \in \mathcal{Z}_{0}^{(n_{\max}-1)} \right\} \\
&=  \left\{ z \;|\; \overline{A}_{p-1} z \in \mathcal{Z}_{0}^{(n_{\max})} \right\},
\end{split}
\end{equation}
which implies that
\begin{equation} \label{equ4_n_max_1}
\begin{split}
\forall z \in \mathcal{Z}_{p-1}^{(n_{\max})}, \; \overline{A}_{p-1} z \in \mathcal{Z}_{0}^{(n_{\max})}.
\end{split}
\end{equation}
Similarly, step 6 in Algorithm~\ref{algorithm_polytopic} implies that
\begin{equation} \label{equ4_n_max_2}
\begin{split}
\forall z \in \mathcal{Z}_{j}^{(n_{\max})}, \; \overline{A}_{j} z \in \mathcal{Z}_{j+1}^{(n_{\max})},\; \forall j \in \mathbb{I}_{[0,p-2]}.
\end{split}
\end{equation}
Combining \eqref{equ4_n_max_1} and \eqref{equ4_n_max_2}, we derive
\begin{equation} \label{equ4_n_max}
\begin{split}
\forall z \in \mathcal{Z}_{j}^{(n_{\max})}, \; \overline{A}_{j} z \in \mathcal{Z}_{j+1\;\text{mod}\; p}^{(n_{\max})},\; \forall j \in \mathbb{I}_{[0,p-1]},
\end{split}
\end{equation}
which together with state constraints satisfaction implies that $\mathcal{Z}^\text{inv}$ is a $p$-periodic invariant tube as defined in Definition~\ref{definition_invariant_tube}, which completes the proof.
\hfill $\qed$
\end{pf}

A sufficient condition for the convergence of Algorithm~\ref{algorithm_polytopic} in a finite number of iterations can be established based on \cite[Theorem~V.1~\textit{(iii)}]{mirceaTAC} applied to a lifted system that admits the lifted limit cycle as a fixed point. The full derivations are however beyond the scope of this paper and will be considered in a separate note focusing on computation of polytopic $p$-periodic invariant tubes for switched affine systems. 

Once an ellipsoidal or polytopic $p$-periodic invariant tube $\mathcal{Z}^{inv}$ is computed, we can simply construct a $p$-periodic invariant tube $\mathcal{X}^{inv}$ in the system state space as shown in Lemma~\ref{lem_Zset}. This can be done analytically for ellipsoidal sets and by computing Minkowski additions with the limit cycle states for polytopic sets.

\subsection{Polytopic outer approximation of $\mathbb{X}_f(N)$}
Consider the proposed FCS-MPC problem \eqref{equ3_FCSMPC_limit_cycle_problem} with a set of $p$ terminal sets $\mathcal{X}_i$.
The set of feasible states $\mathbb{X}_f(N)$ of the FCS-MPC problem \eqref{equ3_FCSMPC_limit_cycle_problem} can be determined by iteratively calculating the one-step controllable set \eqref{equ4_one_step_set} over $N$ iterations, which is defined as follows:
\begin{equation} \label{equ4_one_step_set}
\begin{split}
&\mathbb{X}_f(i+1) := \left\{ x \in \mathbb{X} \;|\; \exists u \in \mathbb{U}: A(u)x+b(u) \in \mathbb{X}_f(i) \right\},
\end{split}
\end{equation}
where $\mathbb{X}_f(0)$ is represented as the union of the $p$ terminal sets of the FCS-MPC problem \eqref{equ3_FCSMPC_limit_cycle_problem}, i.e., 
$\mathbb{X}_f(0) := \bigcup \mathcal{X}^\text{inv}$. Since $\mathbb{U}$ is a \textit{finite control set}, the set of feasible states $\mathbb{X}_f(N)$ is not necessarily convex.

Consider a $p$-periodic invariant tube with polytopic sets, as computed in Algorithm~\ref{algorithm_polytopic}.
For the switched affine system \eqref{equ2_switched_model_dis}, it is estimated that at most a total number of $p\times N_s^N$ enumerations are required to calculate the one step controllable set over a convex polytopic set. Consequently, computing the set of feasible states for the proposed FCS-MPC problem is not scaling well with $N$ and demands substantial storage memory.

\begin{algorithm}[h] 
\caption{Polytopic outer approximation of $\mathbb{X}_f(N)$} \label{algorithm_domain}
\hspace*{\algorithmicindent} \textbf{Initialization:} $\widetilde{\mathbb{X}}_f(0) = \text{conv} \{ \bigcup \mathcal{X}^\text{inv} \}$
\begin{algorithmic}[1]
\For {$i=1:1:N$}
    \For {$j=1:1:N_s$}
    \State $u_j \leftarrow \mathbb{U}\{j\}$
    \Statex $\mathbb{C}_i\{j\} \leftarrow \left\{ x \in \mathbb{X} \;|\; A(u_j)x+b(u_j) \in \widetilde{\mathbb{X}}_f(i-1) \right\}$
    \EndFor
    \State $\widetilde{\mathbb{X}}_f(i) = \text{conv} \{ \bigcup_{j=1}^{N_s} \mathbb{C}_i\{j\}  \}$
\EndFor
\end{algorithmic}
\end{algorithm}
In practice, if the impact of non-convexity is negligible when calculating the set of feasible states, we can instead compute a convex outer approximation of the set of feasible states $\widetilde{\mathbb{X}}_f(N)$ as in Algorithm~\ref{algorithm_domain}, which is tractable. 
Initially, $\widetilde{\mathbb{X}}_f(0)$ is outer approximated as a convex hull of the union of all the terminal sets.
Then for each $u_j \in \mathbb{U}$, a one step controllable set $\mathbb{C}_i\{j\}$ over the set $\widetilde{\mathbb{X}}_f(i-1)$ is computed.
In step 5, $\widetilde{\mathbb{X}}_f(i)$ is formulated by taking the convex hull of the union of $\mathbb{C}_i\{j\}$. This method only requires $N_s\times N$ enumerations to calculate the one step controllable set over a convex polytopic set and $N+1$ enumerations for computing the convex hull of a union sets, thereby significantly reducing the computational burden.

\begin{figure*}[h]
	\centering
	\includegraphics[width=0.9\linewidth]{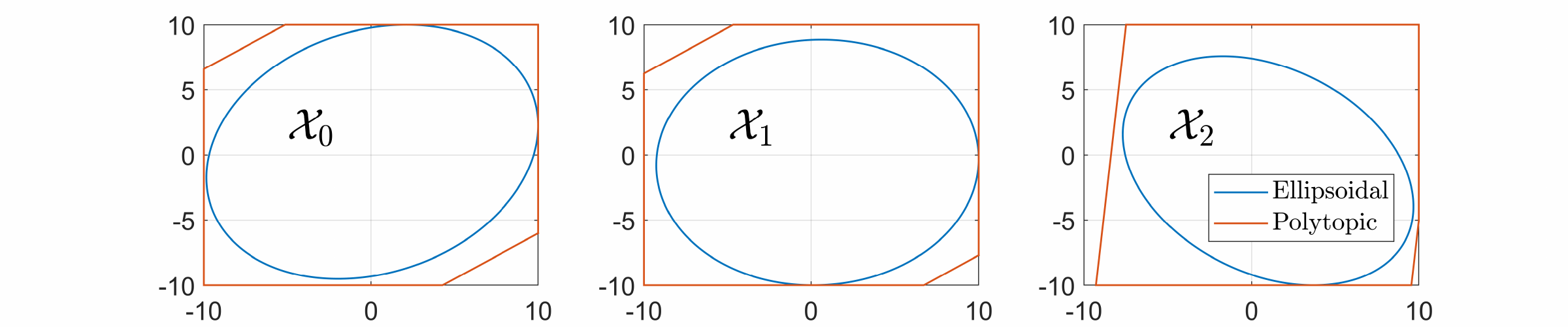}
	\caption{Example 1: Ellipsoidal and polytopic periodic invariant sets.}
	\label{fig1_InvariantSet}
\end{figure*}
\begin{figure*}[h]
	\centering
	\begin{subfigure}[b]{0.4\linewidth}
		\centering
		\includegraphics[width=\linewidth]{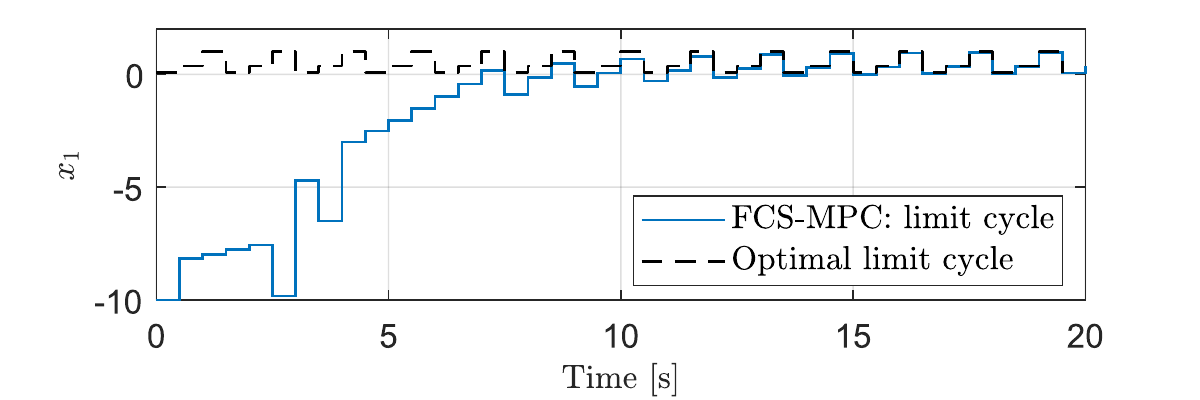}
		\caption{State $x_1$.}
		\label{}
	\end{subfigure}
	\begin{subfigure}[b]{0.4\linewidth}
		\centering
		\includegraphics[width=\linewidth]{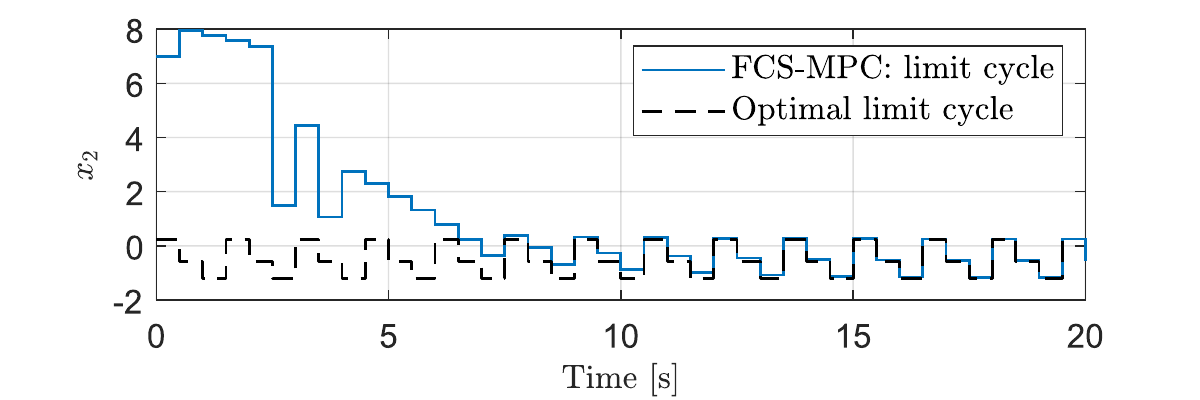}
		\caption{State $x_2$.}
		\label{}
	\end{subfigure}
	\hfill
	\begin{subfigure}[b]{0.4\linewidth}
		\centering
		\includegraphics[width=\linewidth]{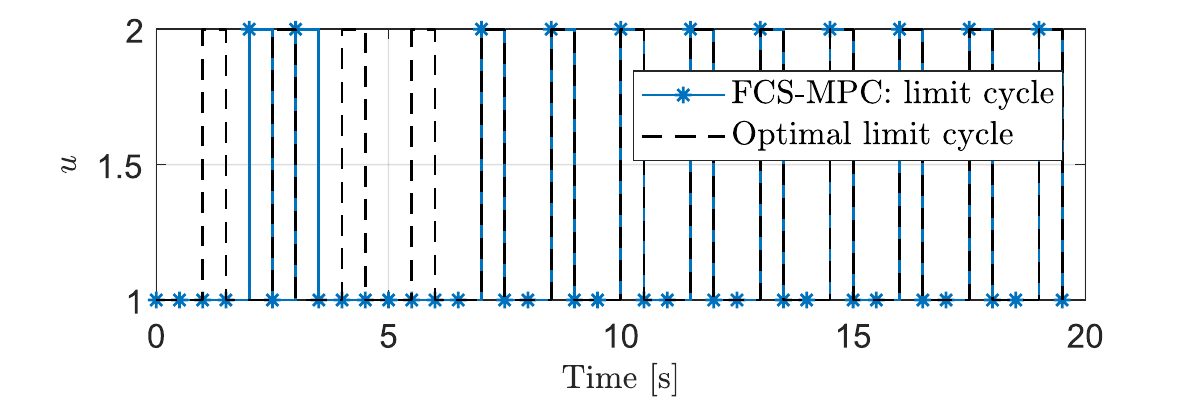}
		\caption{Input $u$.}
		\label{}
	\end{subfigure}
	\begin{subfigure}[b]{0.4\linewidth}
		\centering
		\includegraphics[width=\linewidth]{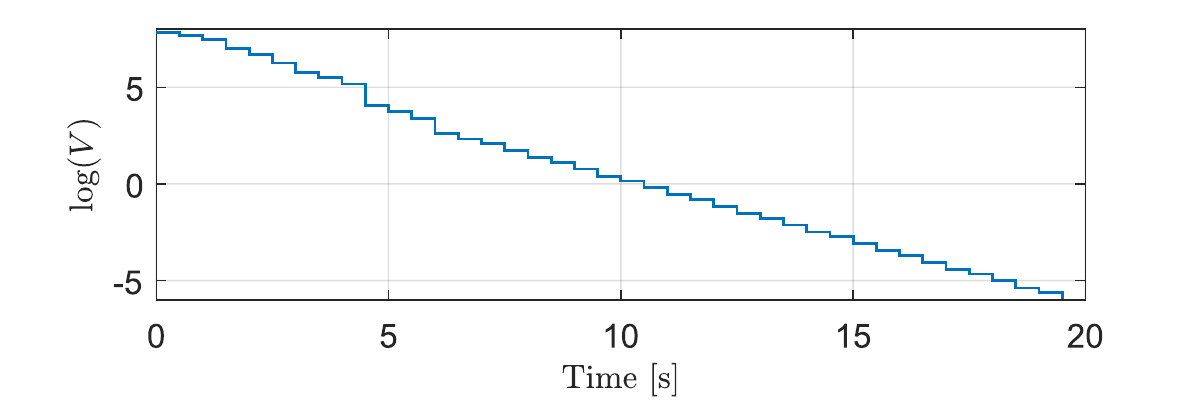}
		\caption{Value function $\log(V)$.}
		\label{}
	\end{subfigure}
	\caption{Example 1: Closed-loop limit cycle FCS-MPC system behavior.}
	\label{fig1_system_behaviour}
\end{figure*}
\begin{figure}[h]
	\centering
	\includegraphics[width=0.7\linewidth]{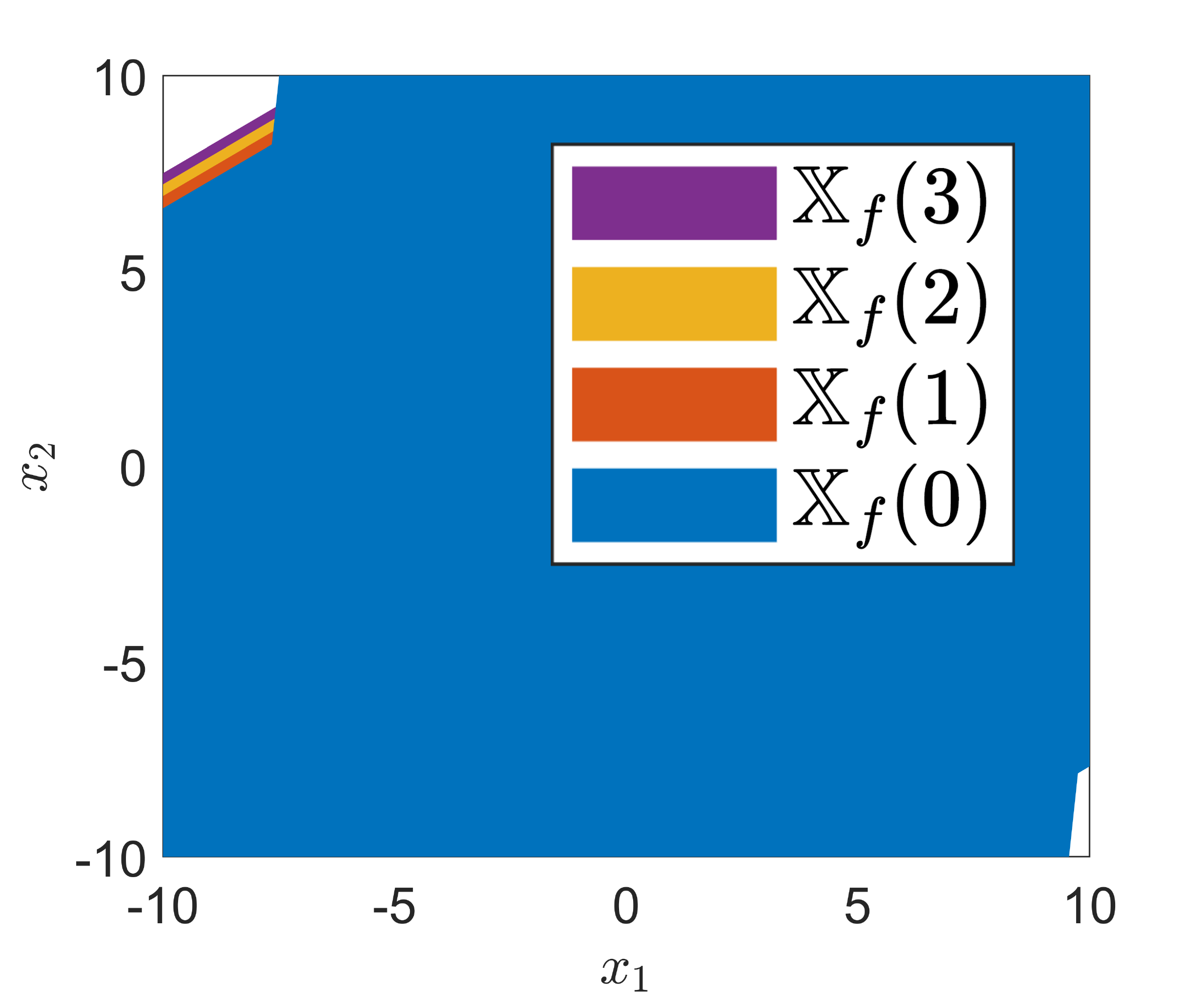}
	\caption{Example 1: Set of feasible states.}
	\label{fig1_feasible}
\end{figure}
\begin{figure}[h]
	\centering
	\includegraphics[width=0.7\linewidth]{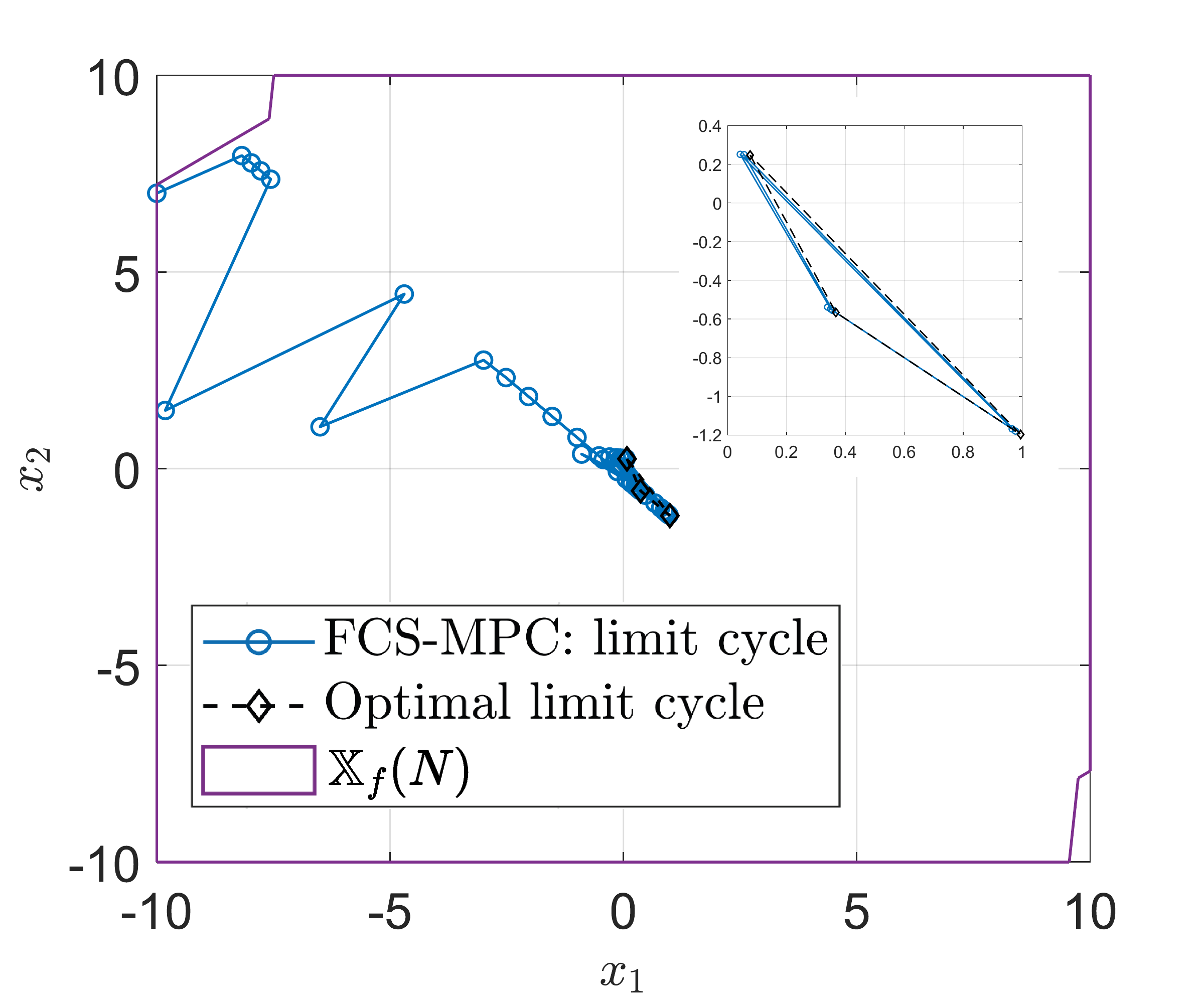}
	\caption{Example 1: State trajectory.}
	\label{fig1_2d_plot}
\end{figure}
\section{Illustrative Examples} \label{Section5}
In this paper, all simulations were conducted using MATLAB 2022b, implemented on a computer equipped with 32GB of RAM and an AMD Ryzen 9 5950X processor, which features 16 cores and 4.9 GHz. The optimal limit cycle is derived by solving the mixed-integer problem \eqref{equ3_limit_cycle_problem} with the cost function $W$ defined in \eqref{equ3_limit_cycle_cost} using the Gurobi solver. The FCS-MPC problem \eqref{equ3_FCSMPC_limit_cycle_problem} is solved online through an exhaustive search approach, executed via parallel computing.

\textbf{Example 1.} Consider the discrete-time switched affine system \cite{egidio2019novel} defined by the matrices
\begin{equation} \label{equ5_model1}
\begin{aligned}
& A_c(u=1)=\begin{pmatrix}-5.8&-5.9\\-4.1&-4.0\end{pmatrix},  &&b_c(u=1)=\begin{pmatrix}0\\-2\end{pmatrix}, \\
& A_c(u=2)=\begin{pmatrix} 0.1&-0.5\\-0.3&-5.0\end{pmatrix},  &&b_c(u=2)=\begin{pmatrix}-2\\2\end{pmatrix}, \\
&C_c =\begin{pmatrix}1&0\\0&1\end{pmatrix}, D_c =\begin{pmatrix}0&0\end{pmatrix},
\end{aligned}
\end{equation}
which represents two unstable subsystems with a sampling time $T_s=0.5$ [s] for discretization. The state and input constraints are defined as
\begin{subequations} \label{}
\begin{align}
\mathbb{X} &= 
\left\{ x | [ -10 \; -10 ]^\top \leq x \leq [10 \; 10 ]^\top \right\},\\
\mathbb{U} &= \mathbb{I}_{[1,2]}.
\end{align}
\end{subequations}
The control target is to steer the system outputs towards the desired reference $\overline{y}=[0 \; 0 ]^\top$. A periodicity $p=3$ is selected to find a corresponding optimal limit cycle as given below:
\begin{subequations} \label{}
\begin{align}
\begin{split}
\overline{X}_{lc}= & \left\{ 
\begin{pmatrix}0.0763\\0.2475\end{pmatrix},\begin{pmatrix}0.3674\\-0.5657\end{pmatrix},\begin{pmatrix}0.9950\\-1.1970\end{pmatrix} \right\},
\end{split}\\
\overline{U}_{lc}= & \left\{ 1,1,2 \right\}.
\end{align}
\end{subequations}
The limit cycle FCS-MPC takes the following setup: $N=4$, $Q=\text{diag}(1,1)$ and $R=0.01$, the set of the terminal weighting matrices $\mathbb{P}_{lc}$ are calculated by \eqref{equ4_P_property}, yielding
\begin{equation} \label{}
\begin{split}
\mathbb{P}_{lc} = \left\{ 
\begin{pmatrix}8.3687&-6.1328\\-6.1328&16.2102\end{pmatrix},\right.\\
\begin{pmatrix}8.8767&-2.9657\\-2.9657&12.1265\end{pmatrix}, 
\left.\begin{pmatrix}14.2377&0.3486\\0.3486&5.6049\end{pmatrix}\right\}.
\end{split}
\end{equation}
Fig. \ref{fig1_InvariantSet} compares the derived ellipsoidal and polytopic terminal sets according to Problem~\ref{problem_ellipsoidal} and Algorithm~\ref{algorithm_polytopic}. Since the computational burden in this setup is not huge, the exact set of feasible states $\mathbb{X}_f(N)$ is revealed in Fig.~\ref{fig1_feasible}.
As demonstrated in Fig. \ref{fig1_system_behaviour} and Fig. \ref{fig1_2d_plot}, the initial state is given as $x(0)=[ -10 \; 7 ]^\top$. Under the designed limit cycle FCS-MPC, all states and the input converge to the optimal limit cycles and satisfy constraints at all times, despite  unstable subsystems.

\begin{figure}[h] 
\centering
\includegraphics[width=0.9\linewidth]{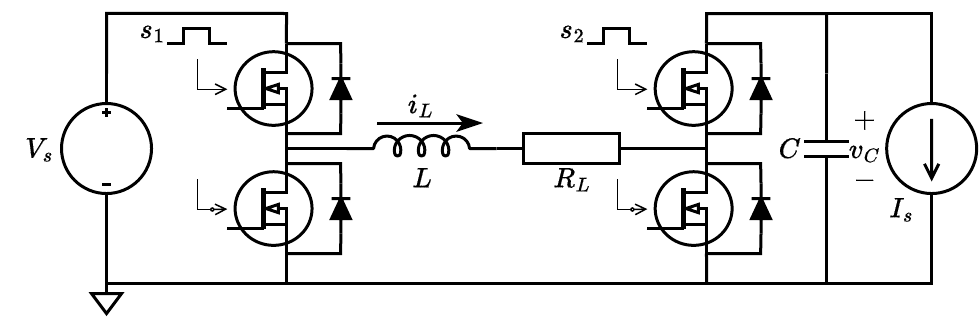}
\caption{Schematic of the Buck Boost converter.}
\label{fig_BB}
\end{figure}
\begin{figure*}[h]
	\centering
	\includegraphics[width=0.8\linewidth]{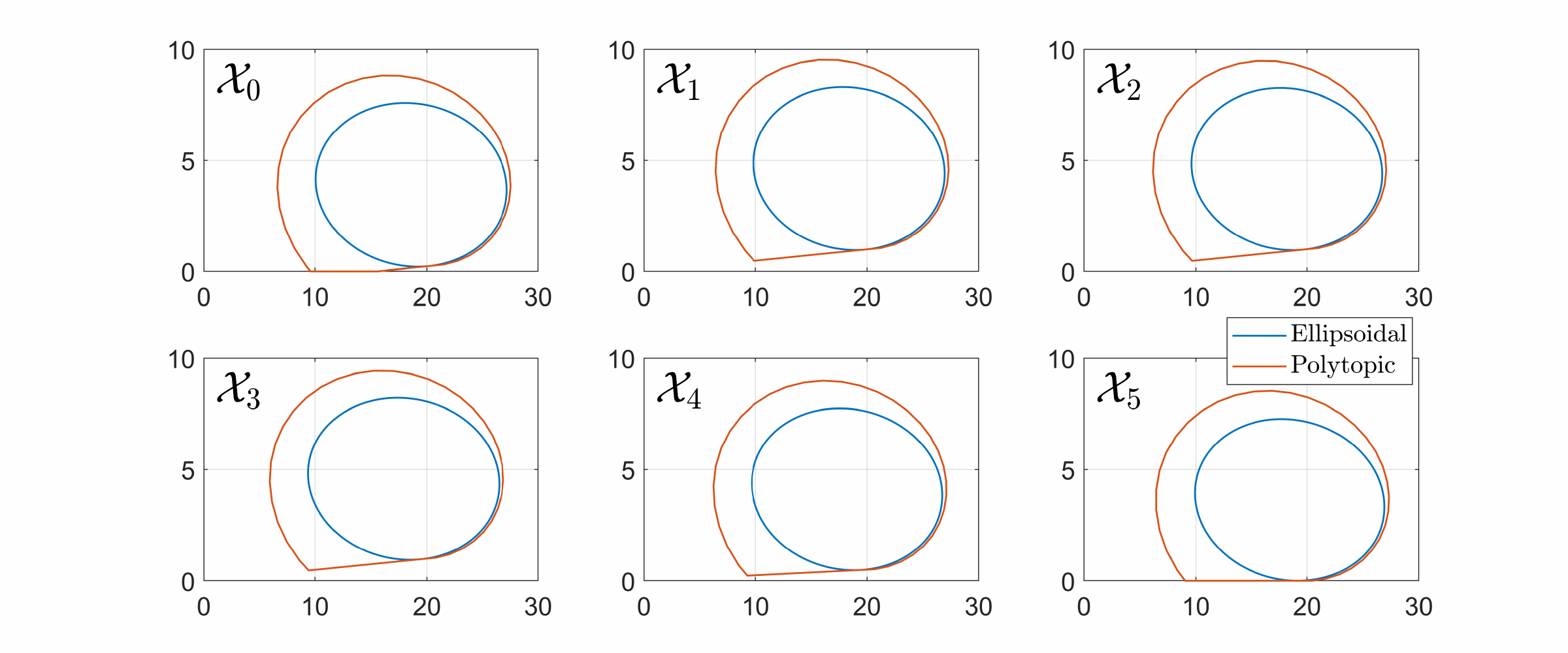}
	\caption{Example 2: Ellipsoidal and polytopic invariant sets.}
	\label{fig_InvariantSet}
\end{figure*}

\textbf{Example 2.} Fig.\ref{fig_BB} illustrates the considered schematic of a non-inverting buck-boost converter borrowed from \cite{xu2022finite}.
Two independent switches $s_1$, $s_2$ control the operation mode of the circuit commutation, which result in $N_s=2^2=4$ modes (or subsystems) in total.
The converter parameters are chosen as $V_s=30$ [V], $I_s=2$ [A], $R_L=0.2$ [$\Omega$], $L=100$ [$\upmu$H] and $C=22$ [$\upmu$F].
\begin{table}[h] 
\centering
\caption{Values of $s_1$ and $s_2$ for each switch mode $\sigma$}
\label{tab_BB_sigma}
\begin{tabular}{rcccc}
\hline
$\sigma$ & 1 & 2 & 3 & 4 \\ \hline
$s_1$     & 0  & 0 & 1 & 1 \\ 
$s_2$     & 0  & 1 & 0 & 1 \\ \hline
\end{tabular}
\end{table}
This topology can be modelled as a continuous-time switched affine system \eqref{equ2_switched_model_con} with $x=[ v_{C} \;i_{L} ]^\top$, $u=[ s_{1} \; s_{2} ]^\top$, $y=v_C$, $\omega=[ V_{s} \; I_{s} ]^\top$ and the system matrices are given below:
\begin{equation} \label{equ5_BB_model}
\begin{aligned}
&A_c(u(t))  =\begin{pmatrix}0&\frac{s_2(t)}{C}\\-\frac{s_2(t)}{L}&-\frac{R_L}{L}\end{pmatrix}, 
&&B_c(u(t)) =\begin{pmatrix}0&-\frac{1}{C}\\\frac{s_1(t)}{L}&0\end{pmatrix}, \\
&C_c(u(t))  =\begin{pmatrix}1&0\end{pmatrix},
&&D_c(u(t)) =\begin{pmatrix}0&0\end{pmatrix}.
\end{aligned}
\end{equation}
The discrete-time switched affine model is obtained considering a sampling frequency of $f_s=400$ [kHz].
The state and input constraints are given as
\begin{subequations} \label{equ5_constraints}
\begin{align}
\mathbb{X} &= 
\left\{ x | [ 0 \; 0 ]^\top \leq x \leq [50 \; 10 ]^\top \right\},\\
\mathbb{U} &= \mathbb{I}_{[0,1]}^2.
\end{align}
\end{subequations}
For simplicity, the switch mode $\sigma$ is used to represent the switch position $u$, as defined in Table \ref{tab_BB_sigma}. 
\begin{figure*}[h]
	\centering
	\begin{subfigure}[b]{0.40\linewidth}
		\centering
		\includegraphics[width=\linewidth]{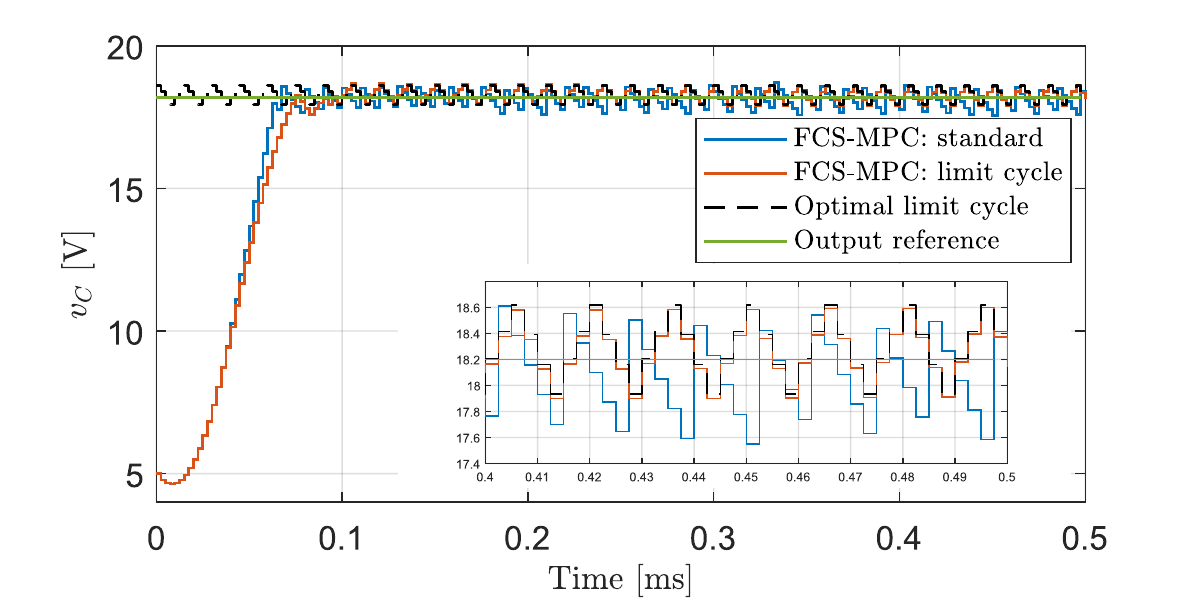}
		\caption{State $x_1$ (equal to the output)}
		\label{fig_x1}
	\end{subfigure}
	\begin{subfigure}[b]{0.40\linewidth}
		\centering
		\includegraphics[width=\linewidth]{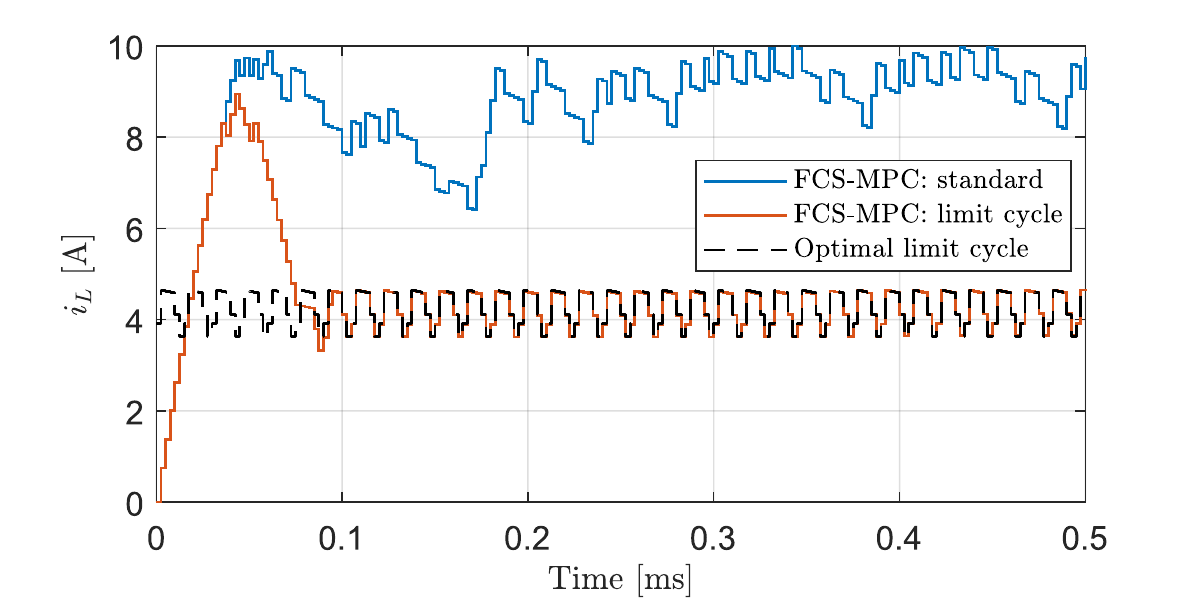}
		\caption{State $x_2$}
		\label{fig_x2}
	\end{subfigure}
	\hfill
	\begin{subfigure}[b]{0.40\linewidth}
		\centering
		\includegraphics[width=\linewidth]{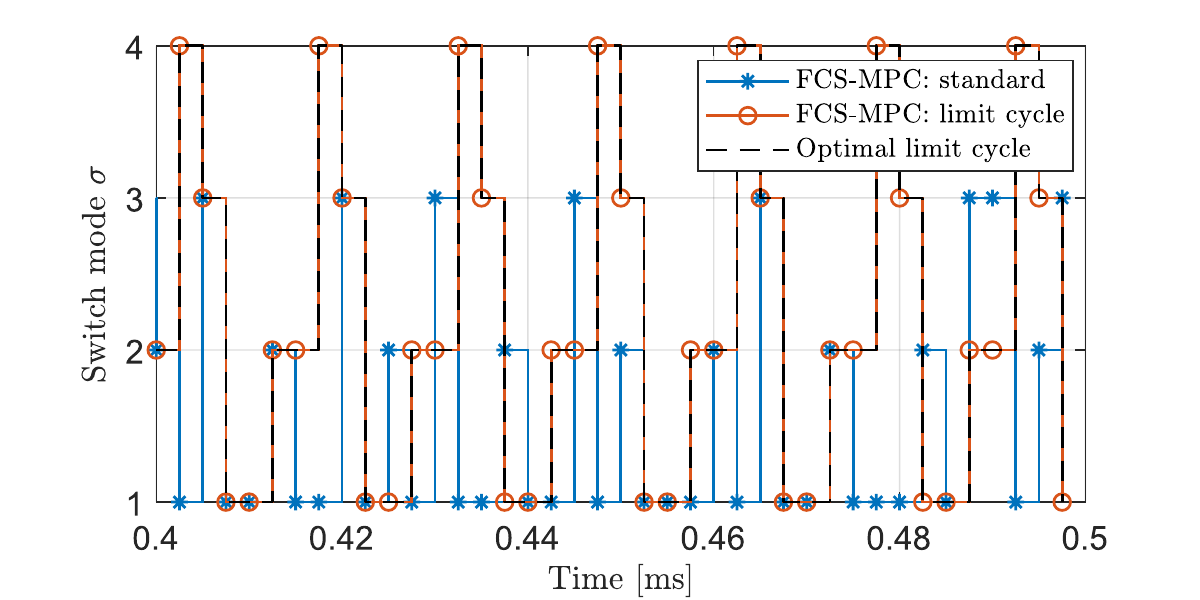}
		\caption{Input $\sigma$ behaviour at steady-state.}
		\label{fig_sigma}
	\end{subfigure}
	\begin{subfigure}[b]{0.40\linewidth}
		\centering
		\includegraphics[width=\linewidth]{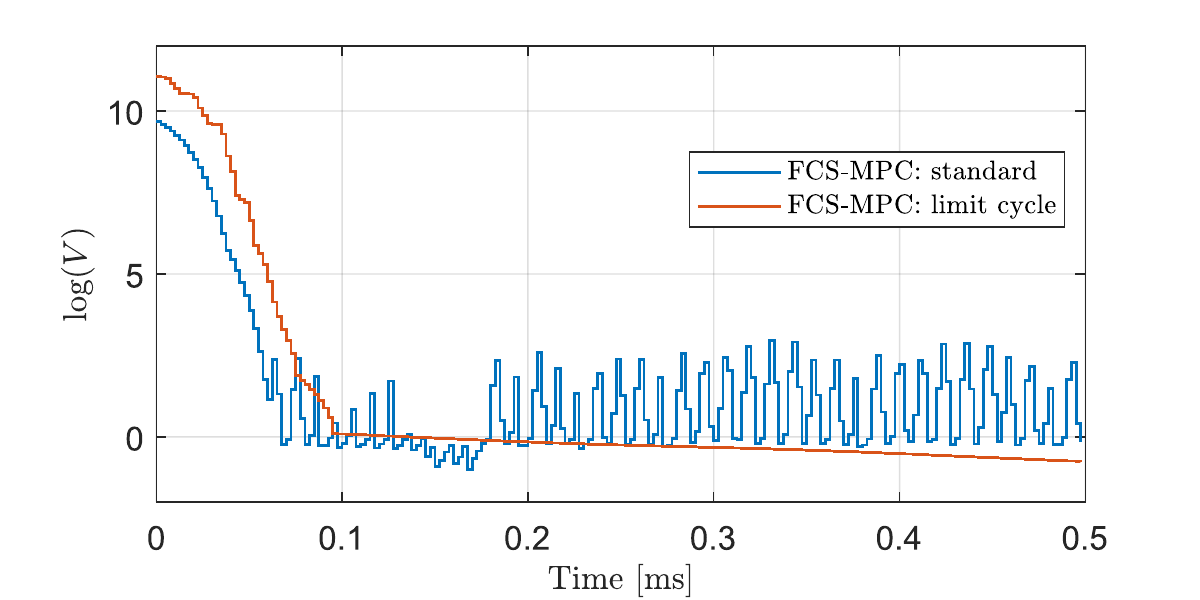}
		\caption{Value function $\log(V)$}
		\label{fig_V}
	\end{subfigure}
	\caption{Example 2: Closed-loop system behaviour; standard FCS-MPC versus limit cycle FCS-MPC.}
	\label{fig_system_behaviour}
\end{figure*}
Considering regulation of the Buck Boost converter towards a constant voltage reference $\overline{y}=v_\text{ref}=18.2$ [V], a standard FCS-MPC is formulated as a benchmark, which has the following cost function with a prediction horizon $N=10$:
\begin{equation} \label{}
\begin{split}
\hat{J}\left( x(k), U_{k}\right) =& \sum_{i=0}^{N-1} \| y_{i|k}-\overline{y} \|^2 + 0.01\| u_{i|k}-u_{i-1|k} \|^2\\
&+ 100 \| y_{N|k}-\overline{y} \|^2.
\end{split}
\end{equation}
For the proposed limit cycle tracking FCS-MPC, the optimal limit cycle in accordance with the reference $\overline{y}$ is derived below by selecting periodicity $p=6$, i.e.,
\begin{subequations} \label{}
\begin{align}
\begin{split}
\overline{X}_{lc}= & \left\{ 
\begin{pmatrix}18.3900\\4.6343\end{pmatrix},\begin{pmatrix}18.1627\\4.6112\end{pmatrix},\begin{pmatrix}17.9355\\4.5882\end{pmatrix}, \right.\\
& \left. \begin{pmatrix}18.2027\\4.1146\end{pmatrix},\begin{pmatrix}18.4159\\3.6374\end{pmatrix},\begin{pmatrix}18.6173\\3.9056\end{pmatrix} \right\},
\end{split}\\
\overline{\sigma}_{lc}= & \left\{ 1,1,2,2,4,3 \right\}.
\end{align}
\end{subequations}
The proposed FCS-MPC takes an equivalent prediction horizon $N=10$.
Given the stage cost $Q=\text{diag}(1,L/C)$ and $R=\text{diag}(0.01,0.01)$, the set of the terminal weighting matrices $\mathbb{P}_{lc}$ are calculated by \eqref{equ4_P_property}, yielding
\begin{equation} \label{}
\begin{split}
\mathbb{P}_{lc} = 10^3 \times & \left\{ 
\begin{pmatrix}0.4290&0.0935\\0.0935&1.8432\end{pmatrix},\begin{pmatrix}0.4266&0.0947\\0.0947&1.8539\end{pmatrix}, \right.\\
&\begin{pmatrix}0.4243&0.0959\\0.0959&1.8648\end{pmatrix},\begin{pmatrix}0.4267&0.0951\\0.0951&1.8540\end{pmatrix}, \\
&\left.
\begin{pmatrix}0.4291&0.0939\\0.0939&1.8433\end{pmatrix},\begin{pmatrix}0.4314&0.0922\\0.0922&1.8326\end{pmatrix}\right\}.
\end{split}
\end{equation}
In addition, the derived ellipsoidal and polytopic terminal sets according to Problem~\ref{problem_ellipsoidal} and Algorithm~\ref{algorithm_polytopic} are plotted in Fig. \ref{fig_InvariantSet}. Fig. \ref{fig_feas} illustrates the polytopic outer approximation of the set of feasible states $\widetilde{\mathbb{X}}_f(N)$ based on Algorithm~\ref{algorithm_domain}.
The control performance of limit cycle FCS-MPC is compared with  standard output reference tracking FCS-MPC .
The initial state is chosen as $x(0)=[ 5 \; 0 ]^\top$. Fig. \ref{fig_x1} and and \ref{fig_x2} present the state time plots and Fig. \ref{fig_sigma} depicts the input (switch mode) behavior at steady-state. Standard and limit cycle FCS-MPC schemes obtain expected tracking behavior for the output $v_C$, with the average tracking errors at steady-state of 0.1466 [V] and 0.0618 [V], respectively.
\begin{figure}[h] 
	\centering
	\includegraphics[width=0.9\linewidth]{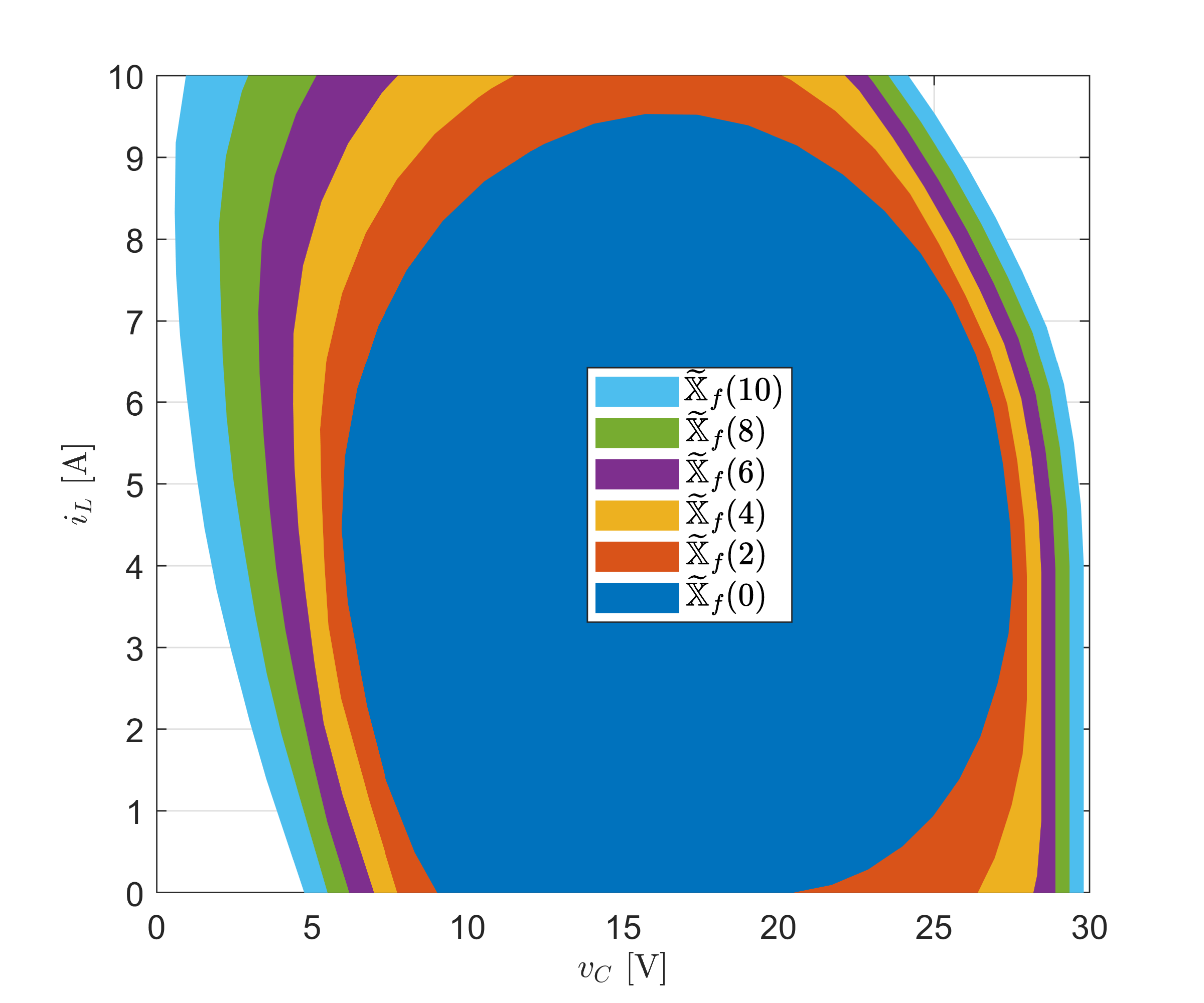}
	\caption{Example 2: Polytopic outer approximation of the set of feasible states.}
	\label{fig_feas}
\end{figure}
\begin{figure}[h] 
\centering
\includegraphics[width=0.9\linewidth]{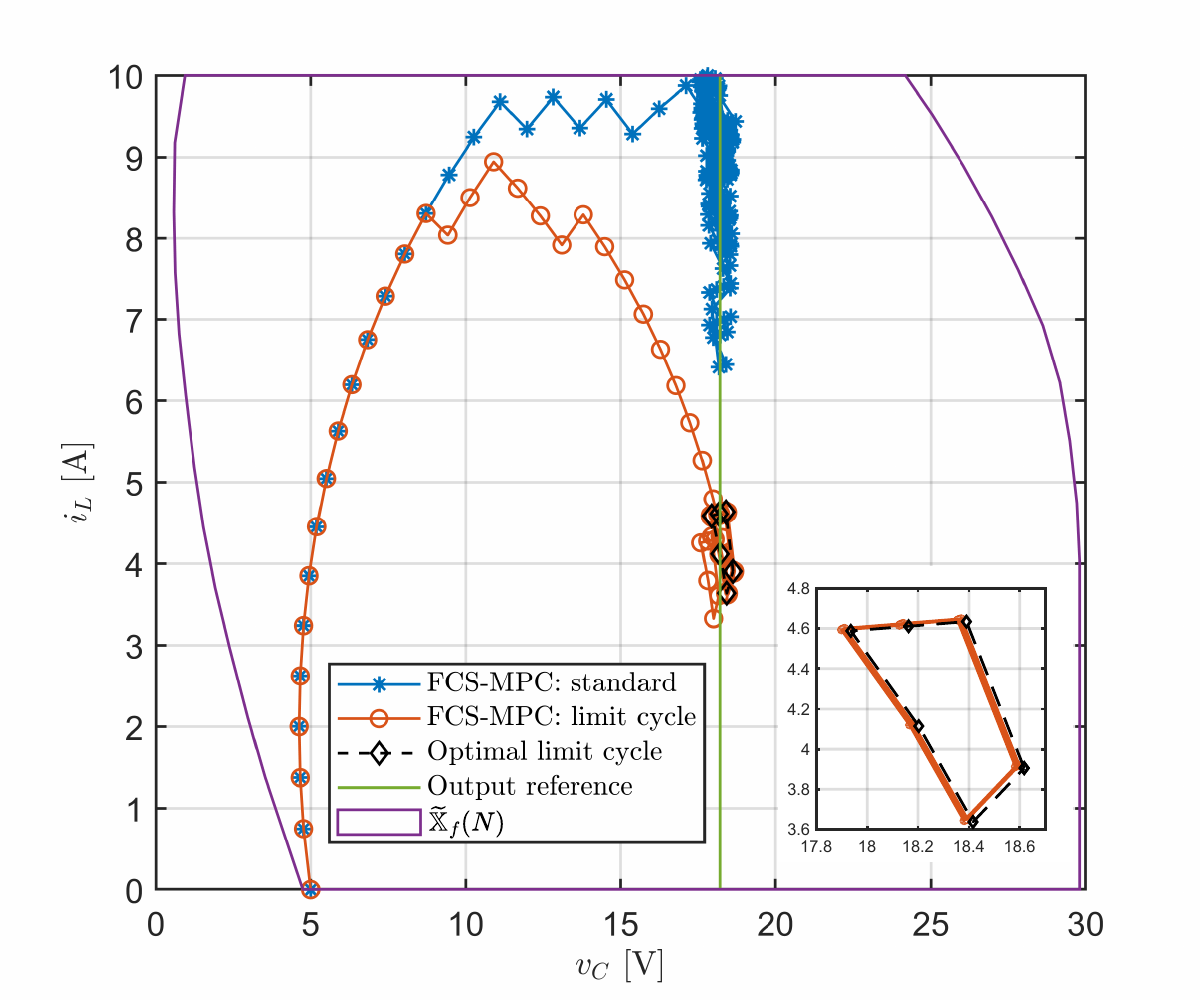}
\caption{Example 2: State trajectory of standard FCS-MPC versus limit cycle FCS-MPC.}
\label{fig_traj}
\end{figure}

The state and input of the proposed FCS-MPC converge to the exact limit cycle as guaranteed. However, at steady-state, the standard output reference tracking FCS-MPC results in a twice higher current value ($i_L$) and no intuitively repetitive performance can be observed, which gives varying switching frequency and is not preferred for power electronics applications. 
Fig.~\ref{fig_V} demonstrates the monotonicity of the optimal value function of the designed limit cycle tracking FCS-MPC.
Fig.~\ref{fig_traj} compares the evolution of the closed-loop state trajectories for standard FCS-MPC versus the developed limit cycle FCS-MPC algorithm.

\section{Conclusions} \label{Section6}
In this paper, a limit cycle FCS-MPC scheme was designed for a class of switched affine systems. We identified a set of assumptions on terminal costs and sets such that the proposed FCS-MPC scheme has guarantees in terms of recursive feasibility and asymptotic stability of the tracking error with respect to a predetermined limit cycle. Approaches for calculating the terminal ingredients and estimating a polytopic outer approximation of the feasible set of states were provided.
The effectiveness of the developed limit cycle FCS-MPC was illustrated on two different examples, one from the literature and a power electronics benchmark converter.

Future work includes the design of robust limit cycle FCS-MPC schemes using, for example, results from \cite{serieye2023attractors} on stabilization of uncertain switched systems, and implementation in real-life power converters.

\begin{ack}  
The authors would like to thank MSc. Sander Damsma for his contributions to the preliminary results on limit cycle FCS-MPC cost function design for linear systems with a finite control set \cite{xu2022steady}, carried out during his MSc graduation project.
\end{ack}

\bibliographystyle{elsarticle-num}
	\bibliography{Section_Bib}

\end{document}